\documentclass[preprint,3p,times]{elsarticle}

\usepackage{lineno,url}
\usepackage{subcaption}
\modulolinenumbers[5]
\usepackage{amsmath,amssymb}
\usepackage{xcolor}
\usepackage[parfill]{parskip}
\usepackage{graphicx}
\usepackage{float}
\usepackage{hyperref}
\usepackage{tikz}
\usepackage{orcidlink}

\newcounter{bla}

\journal{Computer Physics Communications}

\begin{document}

\begin{frontmatter}

%%%%%%%%%%%%%%%%%%%%%%%%%%%%%%%%%%%%%%%%%%%%%%%%%
%%%%%%%%%%%%%%%%%%% Title %%%%%%%%%%%%%%%%%%%%%%%
%%%%%%%%%%%%%%%%%%%%%%%%%%%%%%%%%%%%%%%%%%%%%%%%%

\title{3-Dimensional Adaptive Unstructured Tessellated Look-up Tables for the Approximation of Compton Form Factors}

%%%%%%%%%%%%%%%%%%%%%%%%%%%%%%%%%%%%%%%%%%%%%%%%%
%%%%%%%%%%%%%%%%%% AUTHORS %%%%%%%%%%%%%%%%%%%%%%
%%%%%%%%%%%%%%%%%%%%%%%%%%%%%%%%%%%%%%%%%%%%%%%%%

%% Title, authors and addresses

%% use the tnoteref command within \title for footnotes;
%% use the tnotetext command for the associated footnote;
%% use the fnref command within \author or \address for footnotes;
%% use the fntext command for the associated footnote;
%% use the corref command within \author for corresponding author footnotes;
%% use the cortext command for the associated footnote;
%% use the ead command for the email address,
%% and the form \ead[url] for the home page:
%%
%% \title{Title\tnoteref{label1}}
%% \tnotetext[label1]{}
%% \author{Name\corref{cor1}\fnref{label2}}
%% \ead{email address}
%% \ead[url]{home page}
%% \fntext[label2]{}
%% \cortext[cor1]{}
%% \address{Address\fnref{label3}}
%% \fntext[label3]{}

%% use optional labels to link authors explicitly to addresses:
%% \author[label1,label2]{<author name>}
%% \address[label1]{<address>}
%% \address[label2]{<address>}
 
\corref{*}
\author[a]{Charles Hyde
\orcidlink{0000-0001-7282-8120}
\corref{co-first-author}
%\cortext{corresponding-author}
}

\author[a]{Mitch Kerver
%\corref{co-first-author}
}
\author[b]{Christos Tsolakis
\orcidlink{0000-0002-0656-9631}
\corref{co-first-author}
}
\author[b]{Polykarpos Thomadakis
%\corref{co-first-author}
} 
\author[b]{Spiros Tsalikis
\orcidlink{0000-0001-5113-7195}
\corref{co-first-author}
}
\author[b]{Kevin Garner
\orcidlink{0000-0003-4138-1017}
}
% On the other hand, both Kevin and Angelos should be acknowledged in the Acknowledgement for what they have done.
\author[b]{Angelos Angelopoulos}

\author[c]{Wirawan Purwanto}
 
%\author[b]{Emmanuel Billias}

% \author[b]{Kevin Garner}
\author[d]{Gagik Gavalian
\orcidlink{0000-0002-6738-5457}
}
\author[d]{Christian Weiss
\orcidlink{0000-0003-0296-5802}
}

\author[a,b]{Nikos Chrisochoides
\corref{corresponding-author}
}

\address[a]{Department of Physics, Old Dominion University, Norfolk, VA, USA}
\address[b]{CRTC, Department of Computer Science, Old Dominion University, Norfolk, VA, USA}
\address[c]{Division of Digital Transformation and Technology, Old Dominion University, Norfolk, VA, USA}
\address[d]{Thomas Jefferson National Accelerator Facility, Newport News, VA, USA}

%\cortext[co-first-author]{Authors contributed equally.}
\cortext[corresponding-author]{Corresponding author. \textit{E-mail address:} nikos@cs.odu.edu}

\begin{abstract}
%% Text of abstract
%\note{- Charles} 
We describe an iterative algorithm to construct an unstructured tessellation of simplices (irregular tetrahedra in 3-dimensions) to approximate an arbitrary function to a desired precision by interpolation. The method is applied to the generation of Compton Form Factors for simulation and analysis of nuclear femtography, as enabled by high energy exclusive processes such as electron-proton scattering producing just an electron, proton, and gamma-ray in the final state. While producing tessellations with only a 1\% mean interpolation error, our results show that the use of such tessellations can significantly decrease the computation time for Monte Carlo event generation by $\sim23$ times for $10^{7}$ events (and using extrapolation, by $\sim955$ times for $10^{10}$ events).

\end{abstract}


\end{frontmatter}

%%%%%%%%%%%%%%%%%%%%%%%%%%%%%%%%%%%%%%%%%%%%%%%%%
%%%%%%%%%%%%%%%%% Section 1 %%%%%%%%%%%%%%%%%%%%%
%%%%%%%%%%%%%%%%%%%%%%%%%%%%%%%%%%%%%%%%%%%%%%%%%

\section{Introduction}

\subsection{Nuclear Femtography }
Quantum Chromodynamics (QCD) is a comprehensive theory of the binding of
atomic nuclei, and the structure of subatomic particles including protons and neutrons.
In the past two decades, new theoretical concepts and  experimental techniques have emerged
that for the first time enable the construction of spatial images 
 of quarks and gluons inside the proton and
other atomic nuclei \cite{ji1997gauge}.
The experimental processes are given the generic name Deep Virtual Exclusive Scattering (DVES).
In the case of electron proton scattering ($ep\to eX$), DVES events have the general topology of
\begin{equation}
e + p \to e' + p' + m,
\label{eq:DVESm}
\end{equation}
where $e',\,p'$ refer to the electron and proton in the final state, and $m$ refers to either a meson
(\textit{e.g.,} $\pi^0,\rho,\phi\ldots$) or high energy gamma-ray ($\gamma$) in the final state.

Past and continuing efforts of DVES experiments
include HERA \cite{Aaron:2007ab,Chekanov:2008vy}, HERMES \cite{Airapetian:2012pg}, COMPASS\cite{Akhunzyanov:2018nut}, and Jefferson Lab
\cite{CLAS:2015bqi,CLAS:2018bgk,Defurne:2017paw,Benali:2020vma,JeffersonLabHallA:2020dhq}
%\cite{Kim:2015pkf,Hyde:2011ke,Dlamini:2020ulg}, 
with a major program expected at the U.S. 
Electron Ion Collider (EIC) \cite{Accardi:2012qut} starting in the early 2030s.
For a given luminosity $\mathcal L$ (measured in $\text{cm}^{-2}\text{sec}^{-1}$) of $ep$ collisions,
the event rate for a particular process such as Eq.~\ref{eq:DVESm}, is given by the product of
luminosity times cross section:
\begin{equation}
\text{Rate} = \mathcal L \sigma
\label{eq:DVES}
\end{equation}
Since there are a minimum of three particles to detect in the final state of a DVES process, the cross section
is differential in at least five independent variables---nine momentum variables reduced by four constraints from energy and
momentum conservation.  The functional behavior of the cross section on these five variables is determined by a convolution
of factors determined from the 
fundamental QCD theory, together with a series of functions, called Generalized Parton Distributions (GPDs). GPDs encode the
as yet unknown spatial distributions of quarks and gluons in the target proton (or other target nucleus). The connection between the
GPDs and the cross section on the one hand, and GPDs and spatial imaging on the other hand, will be described in more
detail in the following.

The complexity of the reactions of Eq.~\ref{eq:DVESm} and the even greater complexity of the detectors required to observe them
necessitates extensive simulations of both the functional behavior of the differential cross sections  and the response of the detector.  This process is computationally intensive and must also be flexible.  It is not yet possible to calculate the
GPDs directly from the fundamental QCD theory, so we must rely on models, based on partial understanding of the theory
together with existing data.  Since these models evolve with advances in the field, the computational framework
must be adaptable.

In this paper, we describe algorithms to adaptively and efficiently create tessellated look-up tables of
a set of functions called Compton Form Factors (CFFs). The CFFs are functions of three kinematic variables (variables
determined event-by-event in a DVES experiment).  The CFFs
are computed from GPDs via numerically challenging improper
integrals.  The CFFs are the natural interface between theory and experiment, as they are universal process-independent functions
(\textit{i.e.,} do not depend upon the final state meson in Eq.~\ref{eq:DVESm}) and the  cross section for any process can be directly calculated from the CFFs without
additional numerical integration.  The purpose of the 
tessellation procedure is as follows:
\begin{itemize}
\item Obtain a look-up table such that interpolation on the table will globally reproduce the input model
to any pre-determined precision (\textit{e.g.,} $1\%$).
\item Obtain the look-up table of minimal size in a minimal number of computations of discrete CFF values.  This minimizes
the computation time and storage requirements for each update
or new GPD model.
\end{itemize}

\subsection{Tessellations and Monte Carlo Event Generation}
In order to either predict the sensitivity of a proposed
DVES measurement, or to analyze the data from such
an experiment, it is essential to simulate the response
of the experimental detector to a large statistical
ensemble of simulated pseudo-events.  For a given 
final state channel of Eq.~\ref{eq:DVESm}, millions of
random (Monte Carlo) pseudo-events are generated to comprehensively sample
the experimental phase-space.
Two basic approaches are followed in these studies,
\begin{itemize}
\item
Events are generated uniformly in the multi-dimension
phase space (\textit{e.g.,} 5-D for DVES) and each
event is given a numerical weight proportional to
the differential cross section.
\item
Events are generated by an appropriate algorithm
such that the distribution of pseudo events is
proportional to the differential cross section.
\end{itemize}

%%%%%%%%%%%%%%%%%%%%%%%%%%%%%%%%%%%%%%%%%%%%%%%%%
%%%%%%%%%%%%%%%%% Section 2 %%%%%%%%%%%%%%%%%%%%%
%%%%%%%%%%%%%%%%%%%%%%%%%%%%%%%%%%%%%%%%%%%%%%%%%

\section{Background}

\subsection{Physics} \label{sec:Physics}
Deeply Virtual Exclusive Scattering provides a means to create 3-d tomographic images of nucleons and atomic nuclei. Parton Distribution Functions (PDFs) from Deep Inelastic Scattering experiments have been able to provide information about the longitudinal momentum distribution of partons, pointlike constituents later confirmed to be the quarks and gluons which make up nucleons\cite{Stella_2011,PhysRev.185.1975}. These functions can be generalized to provide transverse spatial distributions of partons as a function of their longitudinal momentum fraction in nucleons (neutrons or protons). These functions are known as Generalized Parton Distributions (GPDs). GPDs can be accessed through Deeply Virtual Compton Scattering(DVCS), where an electron is scattered off a quark within a nucleon and emits a real photon in the final state.
In the $e+N \to e+N+\gamma$ exclusive photo-production reaction, the cross section is the coherent sum of the Bethe-Heitler and Compton amplitudes, making it impossible to distinguish between the two processes experimentally. However, contributions from the Bethe-Heitler process can be computed from QED as a function of the nucleon Form Factors leaving the DVCS amplitude given by:
\begin{equation}
    |T_{DVCS}|^2= \frac{e^6}{y^2Q^2} \{\Sigma c_{n} cos(n\phi_{\gamma\gamma}) + \Sigma s_{n}sin(n\phi_{\gamma\gamma}) \}
\end{equation}

Where $c_n$ and $s_n$ are coefficients involving the CFFs \cite{Belitsky2002}. $\xi$ and $t$ can be accessed experimentally from the process kinematics. However, x is not accessible, but rather integrated over in experiments. Therefore, the quantities that can be extracted from experiments are not the GPDs themselves, but the principal value integral of the GPDs, known as Compton Form Factors (CFFs).

\begin{align}
 \mathcal H_q(\xi,t;Q^2) &=    \int_{-1}^{1} H_{q}(x,\xi,t;Q^2)\left[\frac{}{x-\xi+i\epsilon}
 \pm \frac{1}{x+\xi-i\epsilon}\right] dx 
 \nonumber \\
 &= \mathcal P\int_{-1}^{1}\left[ \frac{H_{q}(x,\xi,t;Q^2)}{x-\xi}\pm \frac{H_{q}(x,\xi,t;Q^2)}{x+\xi} \right]dx 
 +i\pi\left[ -H_{q}(\xi,\xi,t) \pm H_{q}(-\xi,\xi,t) \right]
\end{align}

\subsection{Monte Carlo Event Generators}
Monte Carlo event generators are of central importance in nuclear physics for bridging theory to experiment. As such, many popular Monte Carlo event generators have been made
for high energy $ep$, $pp$, and $AA$ collisions.
Specific generators for $ep$ collisions include Pythia \cite{Sjostrand:2014zea,Pythia8.3},  Herwig, and FOAM. When using a general purpose Monte Carlo event generator, some specific, user implemented, physics model for each event is required. With more complex and computationally intensive models, it can be beneficial to divide the domain into smaller cells. This concept is the basis of both this paper and the FOAM event generator. Although similar in purpose, this project differs greatly in its implementation and benefits from FOAM. FOAM divides the integration domain into hyperrectangles with an algorithm that uses binary splitting, while this paper uses simplices. The benefits of simplices is further explained in section \ref{ssec:tessellations}. FOAM's purpose is to act as a general purpose event generator, while this project focuses on generating tessellations of data to be used by other event generators. Using FOAM requires making calls to the physics model during the event generation. However, once a tessellation of a function is made, it can easily be stored and shared as a lookup table for use in other event generators without calling the physics model.  

\subsection{Tessellations}
Tessellations constitute an integral part of many computational methods.
They provide a discretization of the problem space which enables the
application of simpler methods to each discrete element instead of 
dealing with the whole complexity of the domain at once.
Applications of tessellations span from engineering problems where the
goal is the simulation of air flow around vehicles \cite{johnson1999advanced} to bio-mechanical 
applications \cite{drakopoulos2014parallel} where the focus is the deformation of brain
tissues during surgery. 

In this work, tessellations are utilized as an 
enabling technology that allows for the creation of 
look-up tables, enhancing the Monte-Carlo simulations
used in Compton Form Factor calculations (see section \ref{sec:monte-carlo-simulation-using-lookup-tables}).

\subsubsection{Types of Tessellations} \label{ssec:tessellations}
\begin{figure}[!ht]
    \centering
    \begin{subfigure}[t]{0.24\linewidth}
        \centering
        \includegraphics[width=\linewidth]{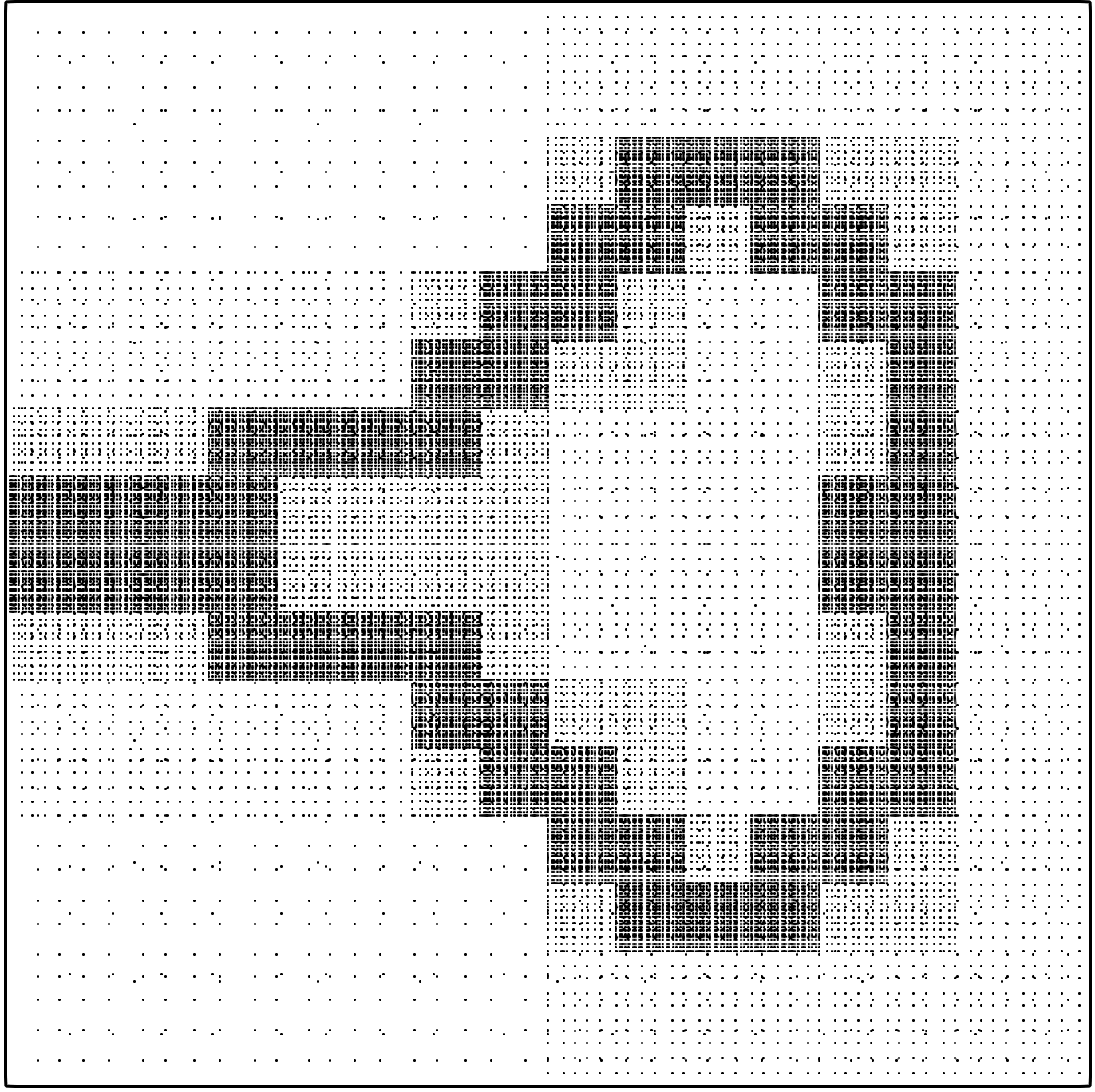}
        \caption{Input Dataset. A uniform tessellation would require 62,000 points in order to capture the smallest features}
        \label{fig:mesh-types-input}
    \end{subfigure}
    \begin{subfigure}[t]{0.24\linewidth}
        \centering
        \includegraphics[width=\linewidth]{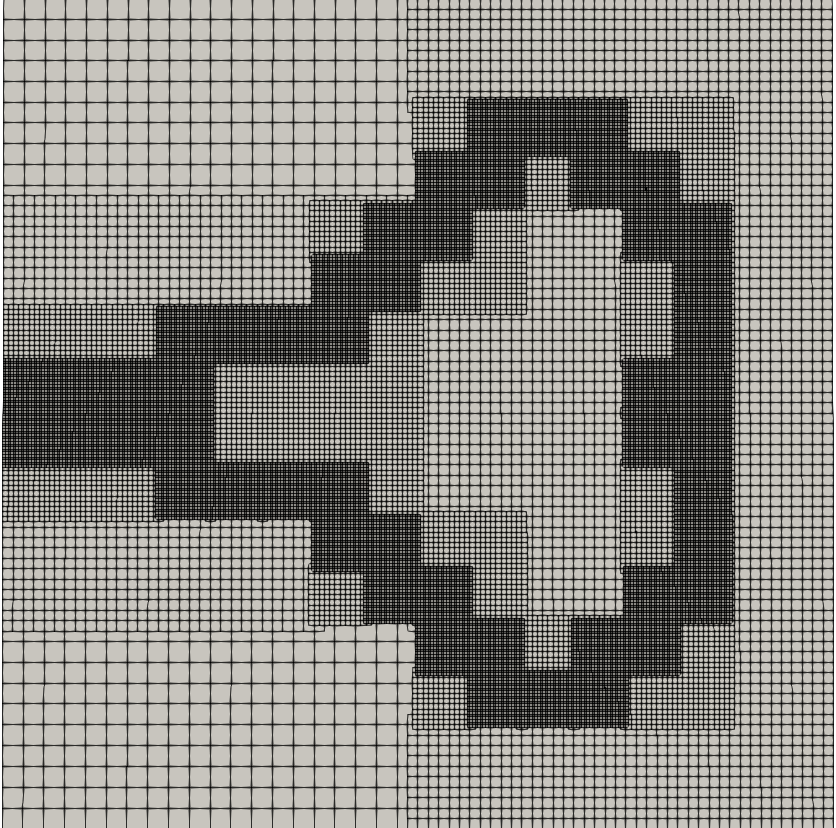}
        \caption{Multiblock tessellation capturing features of the input (59,000 points) }
        \label{fig:mesh-types-amr}
    \end{subfigure}
    \begin{subfigure}[t]{0.24\linewidth}
        \centering
        \includegraphics[width=\linewidth]{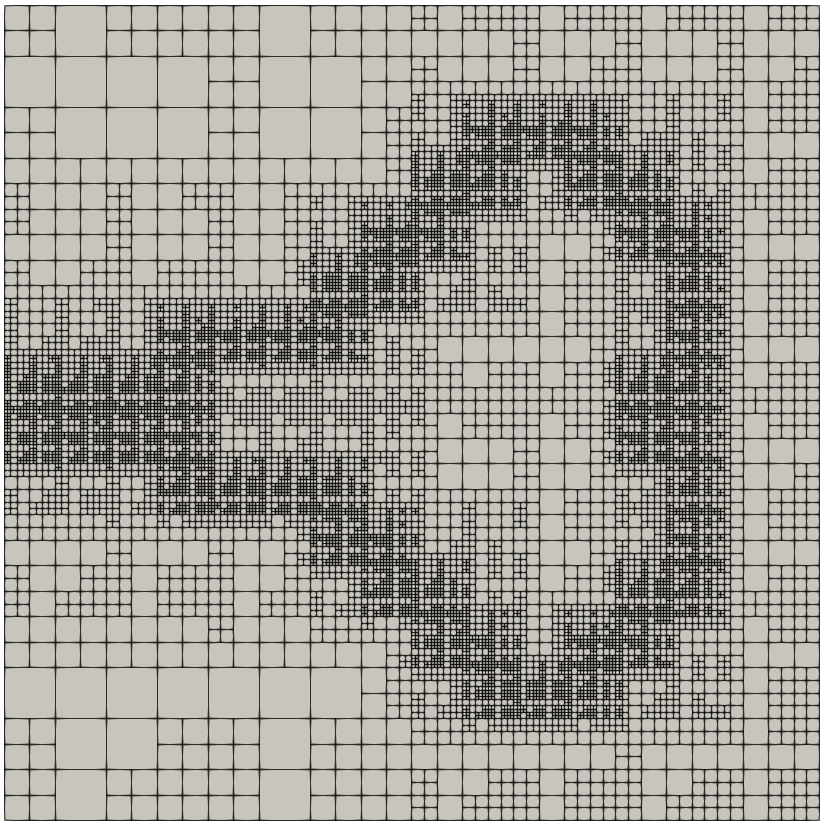}
        \caption{Quad-tree tessellation capturing features of the input (14,500 points) }
        \label{fig:mesh-types-quad}
    \end{subfigure}
    \begin{subfigure}[t]{0.24\linewidth}
        \centering
        \includegraphics[width=\linewidth]{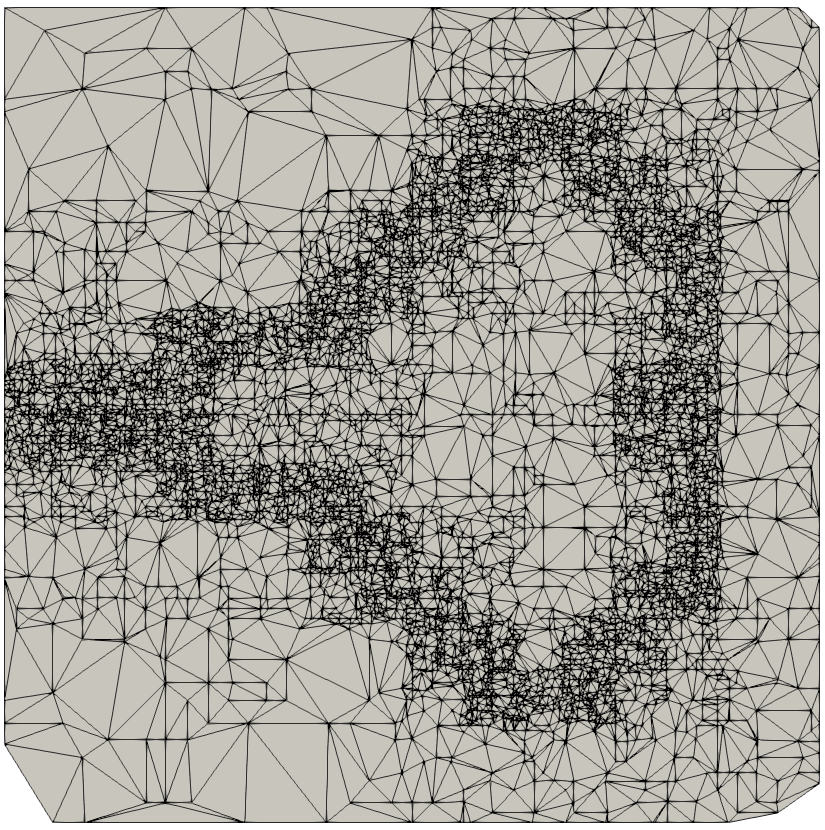}
        \caption{Unstructured simplices (triangles) capturing features of the input (7,500 points) }
        \label{fig:mesh-types-simplex}
    \end{subfigure}
    
    \caption{Different types of tessellations applied on a dataset.}
    \label{fig:mesh-types}
\end{figure}

The plethora of applications that utilize tessellations have different 
requirements and restrictions. Naturally, this has created a wide variety 
of tessellation methods ~\cite{weatherill_handbook_1998}.
Early approaches used structured tessellations 
that arise naturally from the matrix formulation of many problems \cite{weatherill_handbook_1998}. 
The problems of interest rarely have uniform characteristics throughout
the computational space. The ability
to concentrate more points at regions of interest to reduce the discretization error becomes a necessity. Elliptic and hyperbolic methods \cite{chan1992enhancements} provide
a solution to this problem by tweaking the point density at regions of interest. 
Still, they are bound to the initial point distribution and they are limited by 
the complexity of the domain. Multi-block methods (see Figure \ref{fig:mesh-types-amr}) ~\cite{baden_structured_1999} 
provide another way to capture features of the problem by overlaying uniform
tessellations of different density. Octree methods \cite{yerry1984automatic} utilize 
quad-tree (2D) and oct-tree (3D) structures (see Figure \ref{fig:mesh-types-quad}) that can capture 
the features of the problem by using a low count of elements  while at the
same time providing a data structure that favor queries. 
Yet another method  is the use of unstructured simplices (see Figure \ref{fig:mesh-types-simplex})
that not only further reduce the number of points and elements but are
able to scale into higher dimensions at low cost.

\begin{table}[!htp]
\centering
{\scriptsize
\begin{tabular}{l| rrrr | rrrr}
                 & 2-cube & 3-cube & 4-cube & 5-cube & 2-simplex & 3-simplex & 4-simplex & 5-simplex \\
\hline
0-face (vertex)  & 4      & 8      & 16     & 32     & 3         & 4         & 5         & 6         \\
1-face (edge)    & 4      & 12     & 32     & 80     & 3         & 6         & 10        & 15        \\
2-face           & 1      & 6      & 24     & 80     & 1         & 4         & 10        & 20        \\
3-face           &        & 1      & 8      & 40     &           & 1         & 5         & 15        \\
4-face           &        &        & 1      & 10     &           &           & 1         & 16        \\
5-face           &        &        &        & 1      &           &           &           & 1        
\end{tabular}
\caption{Number of lower dimensional faces of n-dimensional cube and simplex elements}
\label{tab:nvertices-growth}
}
\end{table}

Table \ref{tab:nvertices-growth} summarizes the number of lower dimensional elements
contained by each element at different dimensions. Notice that the number of vertices 
scales linearly for simplices while for box-based elements, the growth rate is exponential.
This characteristic is crucial for the approach of this work since it enables the
generalization 
of the method from Compton Form Factor calculations (see section \ref{sec:Physics}) 
to Deeply Virtual Compton Scattering calculations which utilize 5 dimensions and 
Deep Virtual phi-meson Production calculations which can utilize up to 7 dimensions. 

%%%%%%%%%%%%%%%%%%%%%%%%%%%%%%%%%%%%%%%%%%%%%%%%%
%%%%%%%%%%%%%%%%% Section 3 %%%%%%%%%%%%%%%%%%%%%
%%%%%%%%%%%%%%%%%%%%%%%%%%%%%%%%%%%%%%%%%%%%%%%%%

\section{Mesh Generation}\label{sec:mesh-generation}

\subsection{Input}\label{sec:input-image}

\begin{figure}[!ht]
    \centering
    \subfloat[]{
        \includegraphics[width=0.48\linewidth]{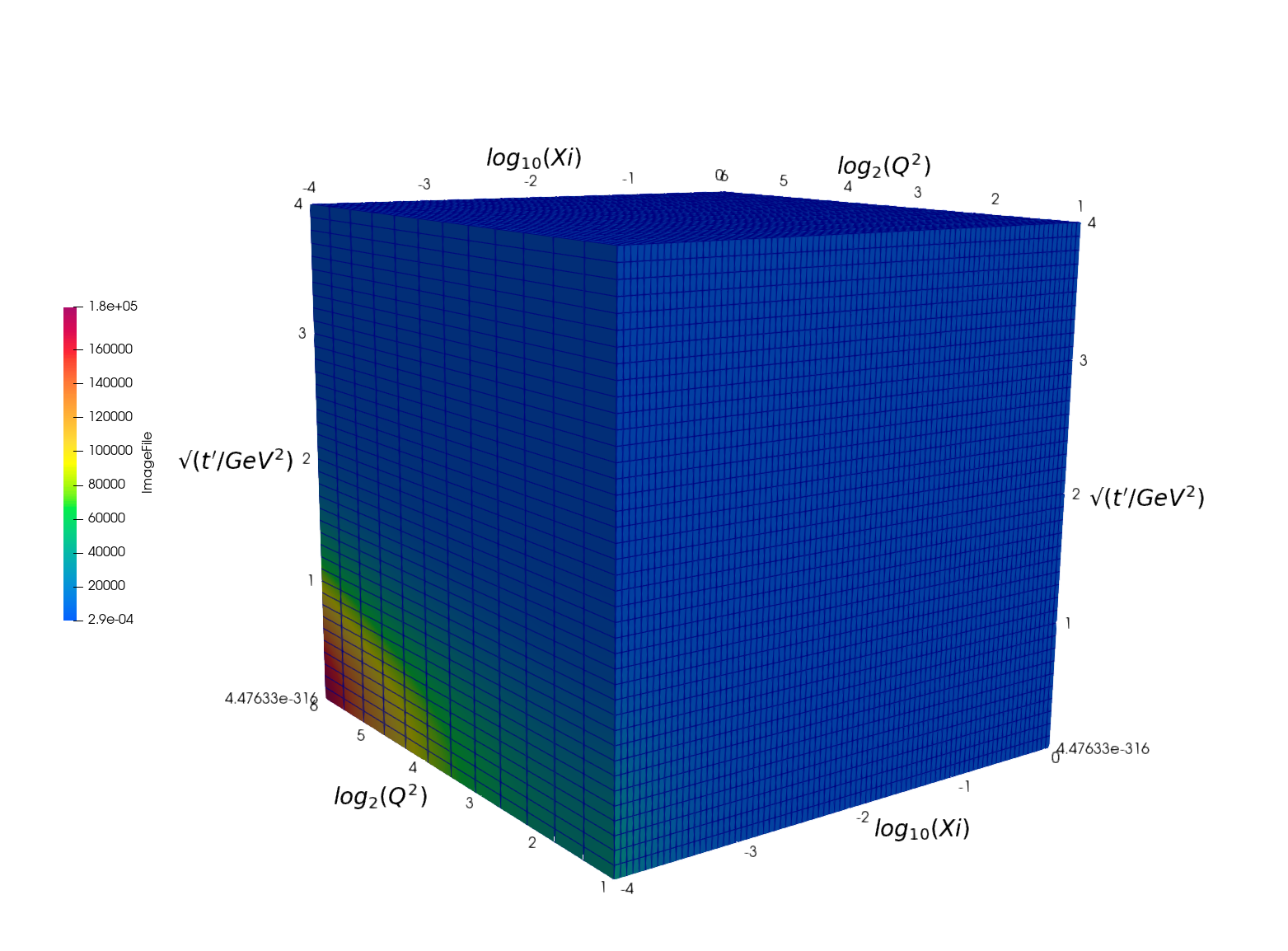}
        \label{fig:structured-imag}
    }
    \subfloat[]{
        \includegraphics[width=0.48\linewidth]{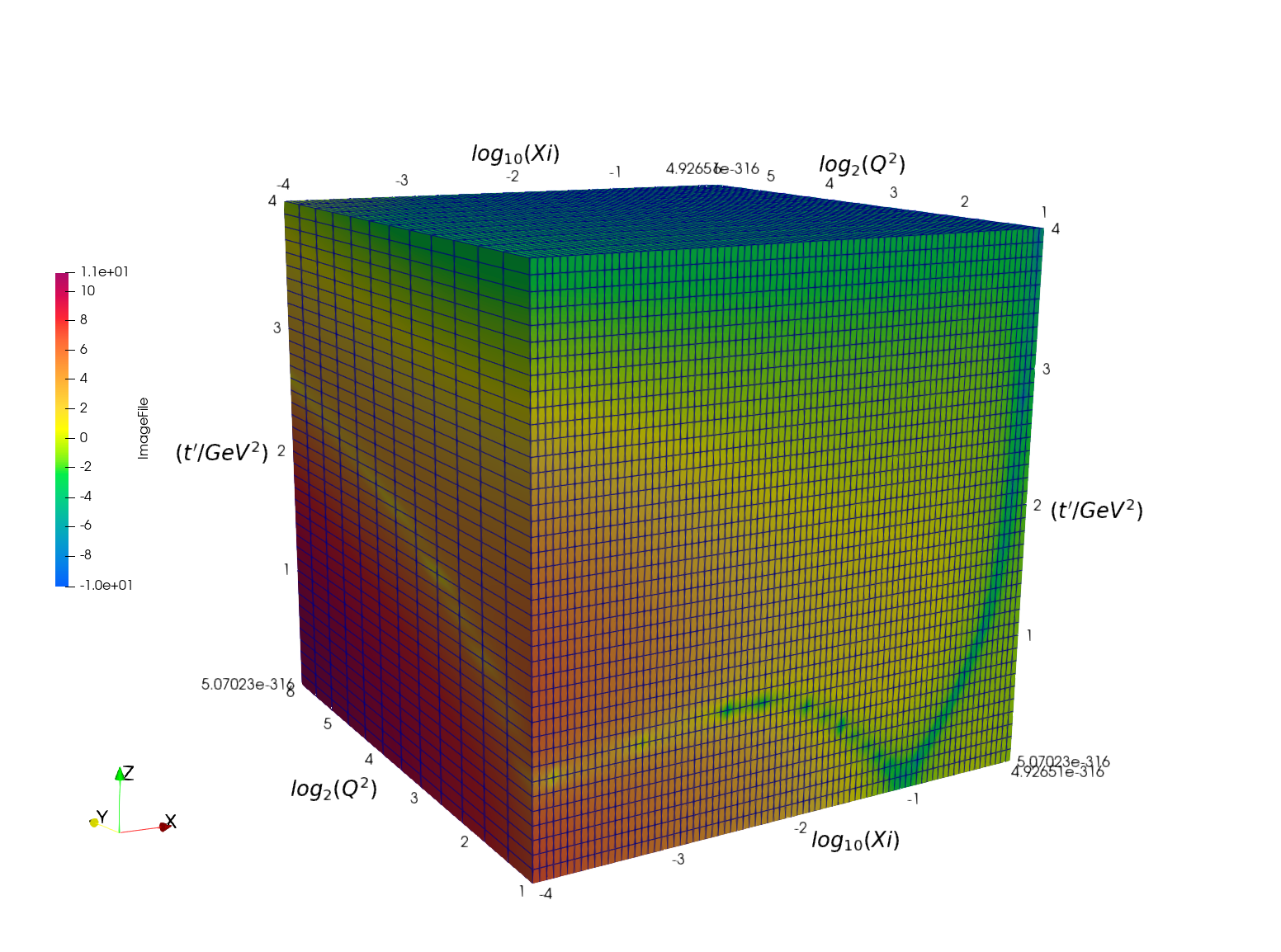}
        \label{fig:structured-real}
    }
    \caption{Cartesian Image of (a) the \textit{Im($H_u$)} CFF and (b) \textit{Re($H_u$)}.  Vertical axis is $\sqrt{t'/0.05 \text{GeV}^2}$. The horizontal axes are $\log_2(Q^2/\text{GeV}^2)$ (y) and $\log_{10}\xi$ (x).}
    \label{fig:PARTONS-NRRD}
\end{figure}

To generate a tessellation of a function, raw data from a model must first be converted to an image-based file format that an image-to-mesh conversion software can understand/process. The file format that we chose is NRRD \cite{NRRD}, because it can support up to 16 dimensions. Once the dimensionality, spacing, and limits of phase space are defined, the model is called and the ITK library \cite{ITK} is used to store the value of the function in an NRRD file as points on a regular Cartesian lattice.  Figure~\ref{fig:PARTONS-NRRD} shows a visualization (generated using Paraview \cite{paraview}) of images of the $Im(H_u)$ and $Re(H_u)$ Compton Form Factor generated from \textit{PARTONS}. The data was generated on a 41x13x81 grid in the kinematic variables [$\sqrt{t'/0.5 GeV^{2}}$, $\log_2(Q^2)$, $\log(X_{\xi})$] so that several orders of magnitude could be covered in $Q^2$ and $X_{\xi}$, while maintaining a near cube hyperrectangle for the ease of visualization. The value of the CFF is stored only at the vertices of the grid, but Paraview's linear interpolation scheme is being used to approximate the function continuously over the entire domain as indicated by the coloring in the visualization.

\subsection{Image-to-mesh conversion software}\label{sec:I2M_software}
The image-to-mesh conversion software utilized is called PODM \cite{PODM}. PODM is a state-of-the-art image-to-mesh conversion software capable of generating in parallel meshes of guaranteed quality, while maintaining fidelity with respect to the boundary of the input image. PODM converts images to tetrahedral meshes using proven sampling techniques. The sampling of the mesh (volume and surface) is controlled by the parameter delta ($\delta$), where $\delta > 0$. A low value of $\delta$ indicates that the sampling is denser. Ideal values of $\delta$ are considered to be multiples of the spacing of the image. To guarantee quality and fidelity, PODM employs a set of rules that split specific elements when needed. In conjunction with these rules, a user-defined rule can be specified as a sizing function, which is also used to split elements. The user-defined rule can be utilized to adapt the mesh to the input data, and it can be loaded as a shared library that PODM will utilize on-demand.

\subsubsection{Portability}\label{sec:portability}
A Docker image containing \textit{PODM} and \textit{PartonsToNRRD} has been created for easy distribution to users. The software operates directly from within a Docker container and all processing happens within it. Data is transferred from a user machine to the Docker container via the use of bind mounts or volumes. Input is read from a user’s computer, processed inside the Docker container, and then the output is written back to the user’s computer. There are many benefits to this approach, including portability, ease-of-use, and isolation (\textit{i.e.,} the isolation offered by Docker technologies, meaning no files are created on the base system). In addition, the usage of Docker ensures the creation of a well-defined runtime environment in which the software has been tested and verified by software developers. Users can simply load the Docker image, spawn a Docker container, and start using the software immediately. The Docker image contains a lightweight distribution of Debian Linux, \textit{PODM}, utilities, and all the libraries required for the software to run. The Docker image can be used on all operating systems that Docker supports, including Windows, Linux, and macOS. More information regarding their repositories and how one can use these methods is described in the appendix.

\subsection{Interpolation Error}\label{sec:interpolation-error}

Tessellations, as discussed in section 2, provide an efficient way to store data with a few number of vertices. During Monte Carlo event generation, this data needs to be quickly accessed at any given point in the domain, which requires an interpolation method to be used between vertices of the tessellation. The linear interpolation method used is described in detail in section \ref{sec:numerical-operations}. Approximating the value of a non-linear function with linear interpolation introduces an error. The method which generates the tessellations and later refines them should aim to minimize the global interpolation error. Therefore, before the different approaches to mesh generation are explored, the method of evaluating interpolation error must be defined.

To evaluate the significance of the interpolation error, sample points are generated at the barycenter of each simplex in the tessellation. At each sample point, the function is evaluated by the GPD model and compared to the linear interpolation value using the relative error, \textit{i.e.,} \textit{relative error} $= \frac{\text{PARTONS value} - \text{linear interpolation value}}{max(\text{PARTONS value}, \text{linear interpolation value})}$ . With functions that have both positive and negative values varying over several orders of magnitude, the relative error can surge up to infinity for values close to zero. For these values, it is more useful to look at their absolute error, \textit{i.e.,} \textit{absolute error} $ = \textit{PARTONS value} - \textit{linear interpolation value}$. The interpolation error of each simplex is then defined as $min(\textit{absolute error}, \textit{relative error})$. This guarantees that most values will use the relative error while the absolute error is only used as the function approaches zero.

\subsection{Mesh generation approaches}

In addition to minimizing the number of calls to the GPD model, and thus minimizing the number of tessellation vertices, the tessellation should also be generated in a way that minimizes the global interpolation error. Generating a tessellation with a pre-specified interpolation precision and with a minimal computation cost is a difficult task. For this reason, we tried four different approaches of mesh generation. The first method generates uniform images (grids) while the remaining three methods start with the same initial image of 43,173 points as described in section \ref{sec:input-image}.

\subsubsection{Structured Tessellation}

\begin{figure}[!ht]
    \centering
    \subfloat[]{
        \includegraphics[width=0.5\linewidth]{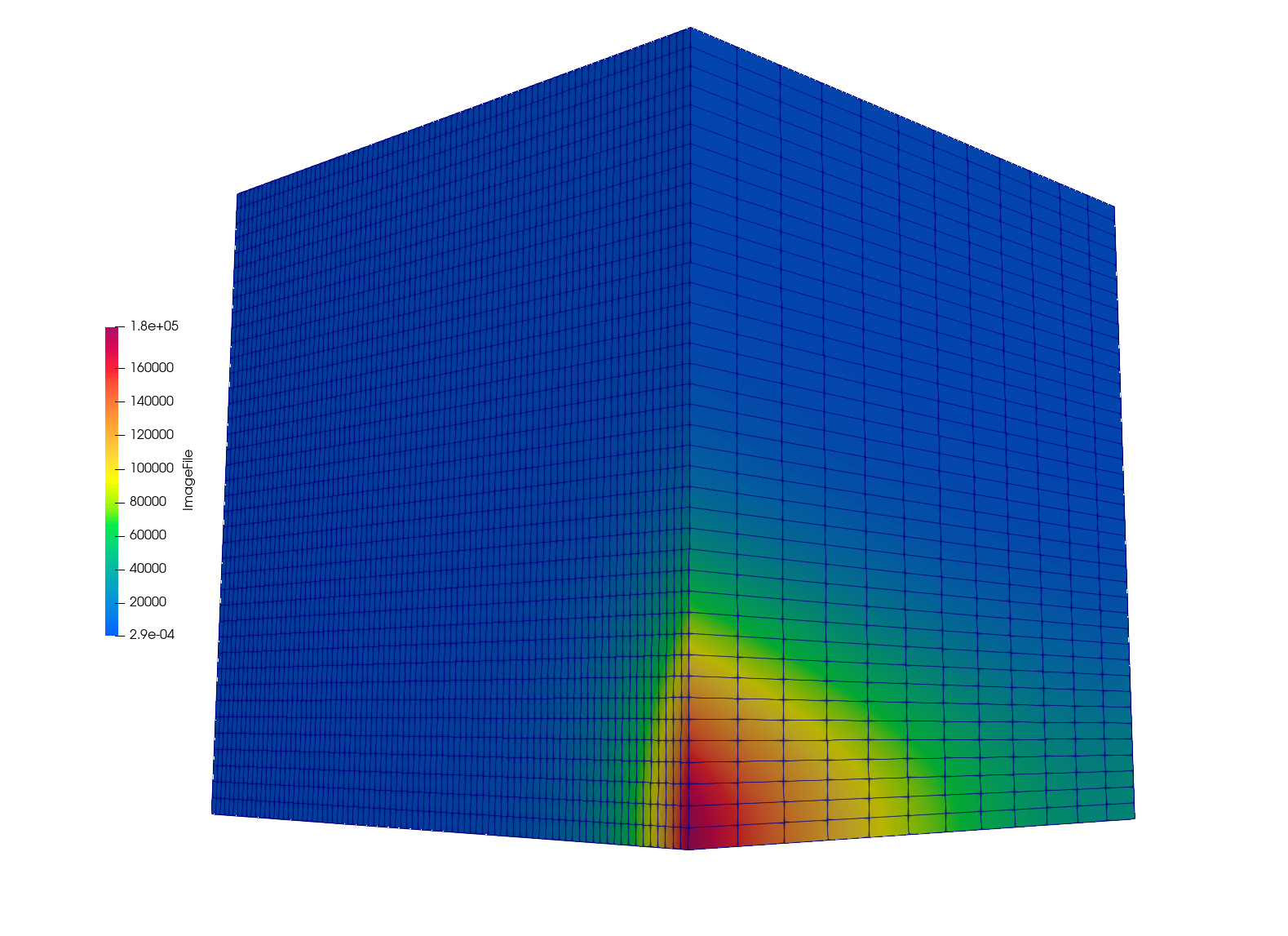}
        \label{fig:zoom-structured-imag}
    }
    \subfloat[]{
        \includegraphics[width=0.5\linewidth]{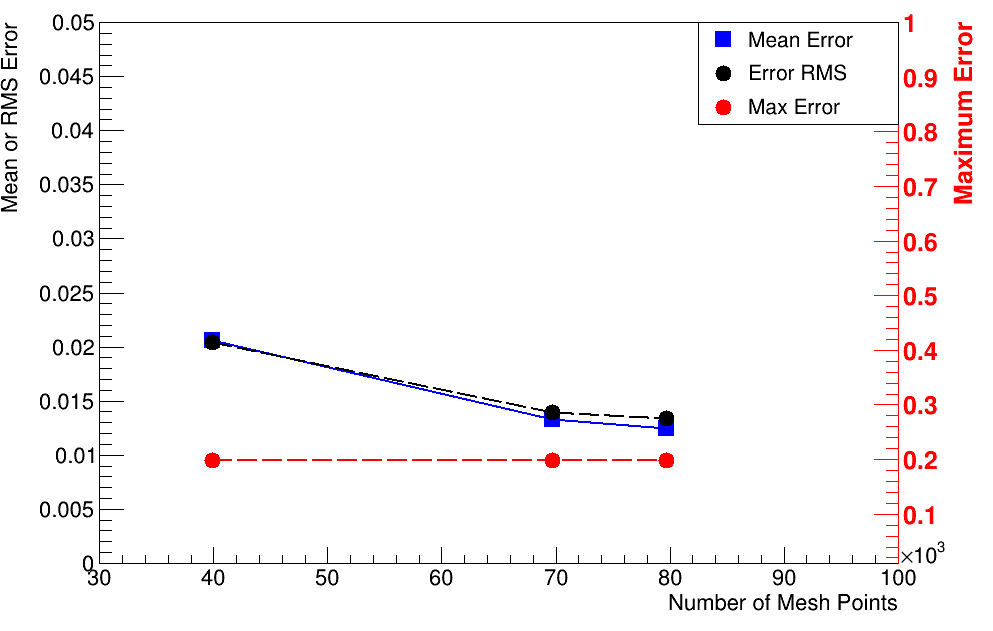}
        \label{fig:meanError_structured}
    }
    \caption{(a) Structured grid of $Im(H_{u})$ CFF and (b) Mean, RMS, and Max interpolation error statistics of structured grids of different sizes}
    \label{fig:structured-imag-meanError}
\end{figure}

In the first approach, structured tessellations (images) are generated with varying numbers of points across each dimension. The more the points across each dimension, the smaller the size of the cubic elements become, thus lowering the interpolation error. Once an initial tessellation is generated, it is checked to see if it meets the desired error threshold. If it does not, the number of points in one or many dimensions is increased. In the case of the $Im(H_u)$ CFF, as seen in Figure \ref{fig:structured-imag-meanError}, we increased the number of points only in the $log_{10}(\xi)$ dimension, because the behavior of the function in $\xi$ is more complex. This loop should continue until the desired precision is met. This approach has several disadvantages: (1) the number of points required across each dimension to meet the desired precision is unknown, (2) a dense \textit{uniform} tessellation has many elements, (3) generating a dense tessellation requires numerous GPD calculations (of which a majority would be wasted on elements that already have a sufficient interpolation error), and (4) for each instance, the GPD model must be called for \textit{all} the points. 

\subsubsection{Unstructured Uniform Tessellation}
\label{sec:unstructured-uniform}

\begin{figure}[!ht]
    \centering
    \subfloat[]{
        \includegraphics[width=0.5\linewidth]{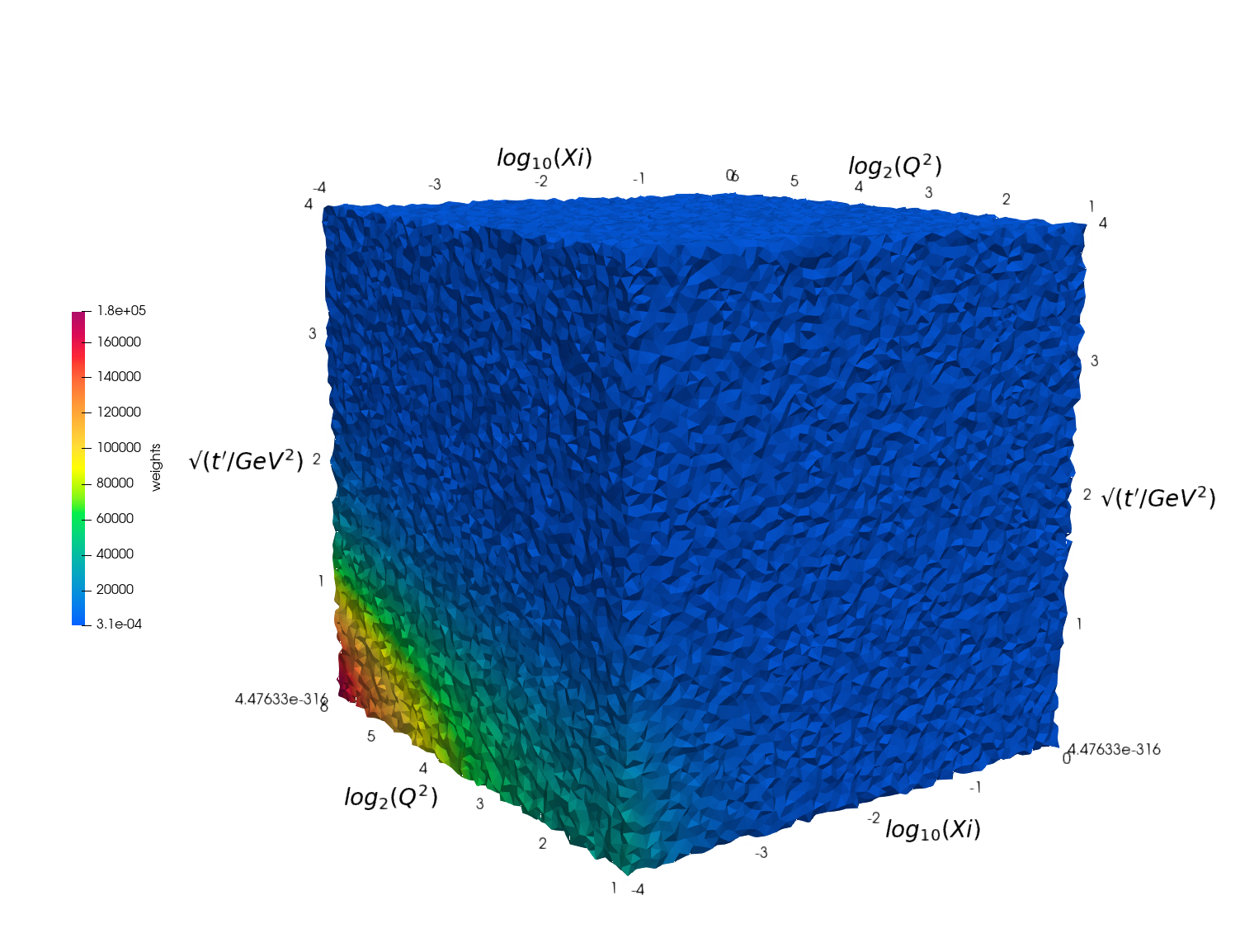}
        \label{fig:uniform-imag}
    }
    \subfloat[]{
        \includegraphics[width=0.5\linewidth]{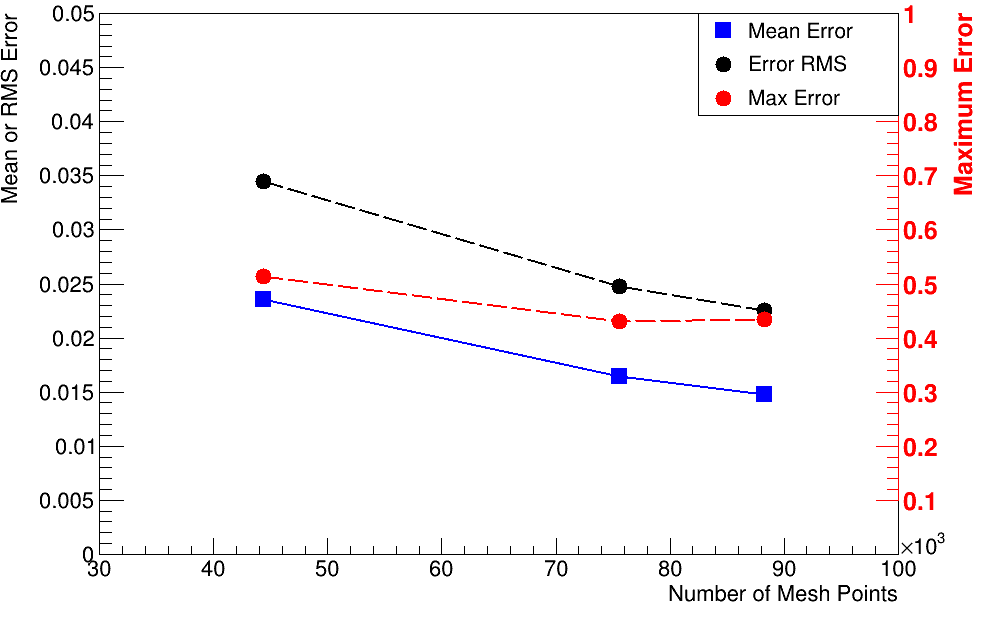}
        \label{fig:meanError_uniform}
    }
    \caption{Uniform tessellation of $Im(H_{u})$ CFF (a) the 1st instance and (b) Mean, RMS, and Max interpolation error statistics of uniform tessellations of different sizes}
    \label{fig:uniform-imag-meanError}
\end{figure}

In the second approach, unstructured uniform tessellations are generated with varying values of $\delta$ to control the density. The denser a uniform tessellation is, the smaller the size of the simplices become, thus lowering the interpolation error. The procedure of generating an unstructured uniform tessellation is similar to the structured approach, and although only $\delta$ needs to be updated for every instance, both have similar disadvantages. This can be seen in Figure \ref{fig:uniform-imag-meanError} which shows the change in the mean and maximum interpolation error as the number of vertices increases. By going from approximately 77,000 vertices to about 90,000 vertices, the mean error decreases slightly. However, the maximum interpolation error remains unchanged at 40\%. This confirms the need for an adaptive tessellation approach, as many points are unnecessarily added in this uniform approach which are wasted in simplices that already have small interpolation error.

\subsubsection{Unstructured Adaptive Tessellation}\label{sec:sizing-functions}

\begin{figure}[!ht]
    \centering
    \subfloat[]{
        \includegraphics[width=0.21\linewidth]{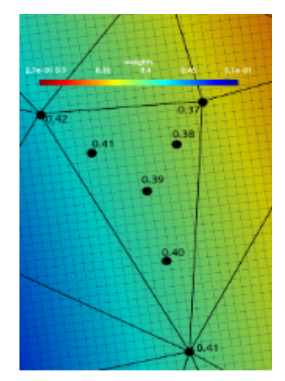}
        \label{fig:sample-points-quering}
    }
    \subfloat[]{
        \includegraphics[width=0.37\linewidth]{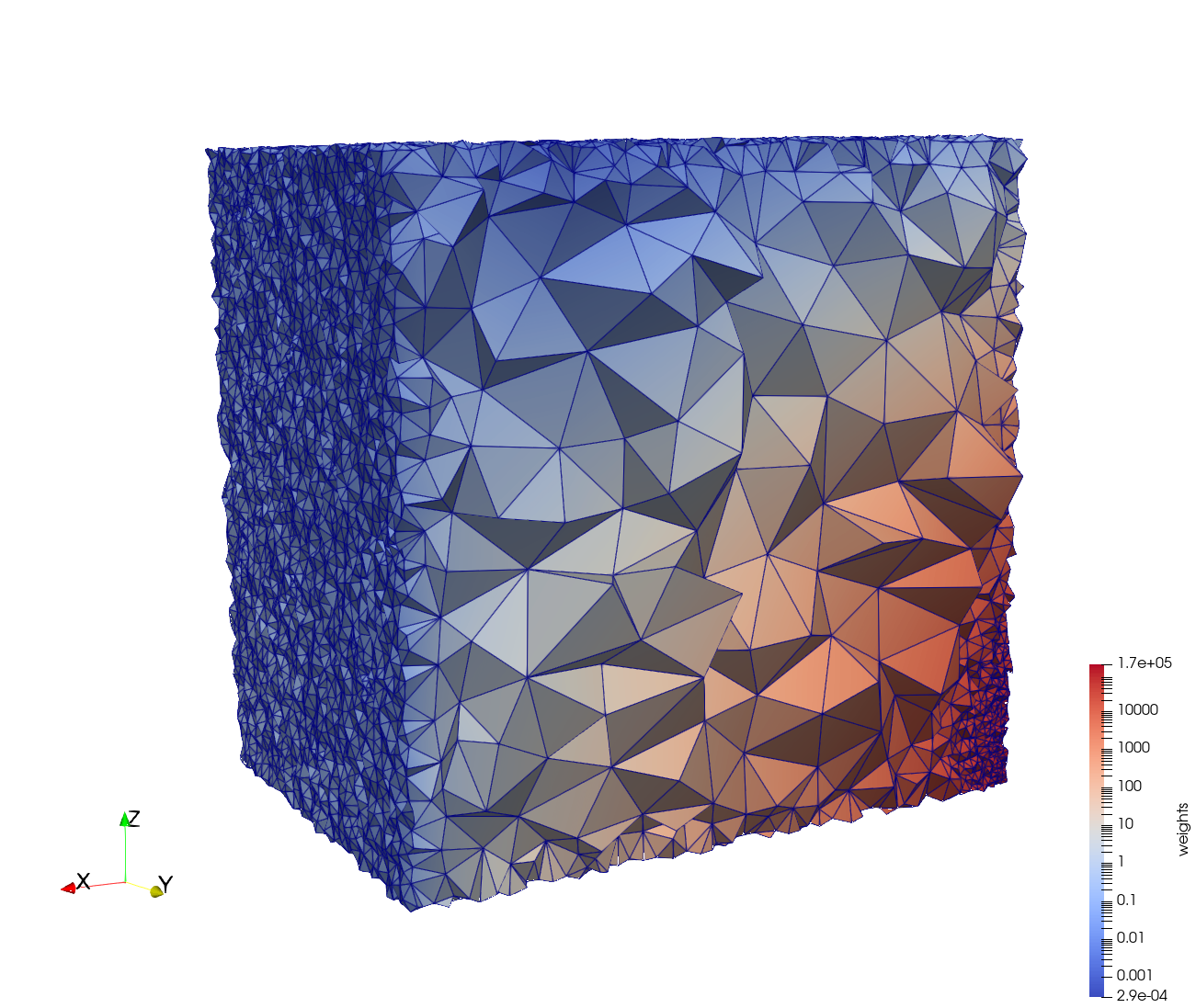}
        \label{fig:adaptative-gradation}
    }
    \subfloat[]{
        \includegraphics[width=0.37\linewidth]{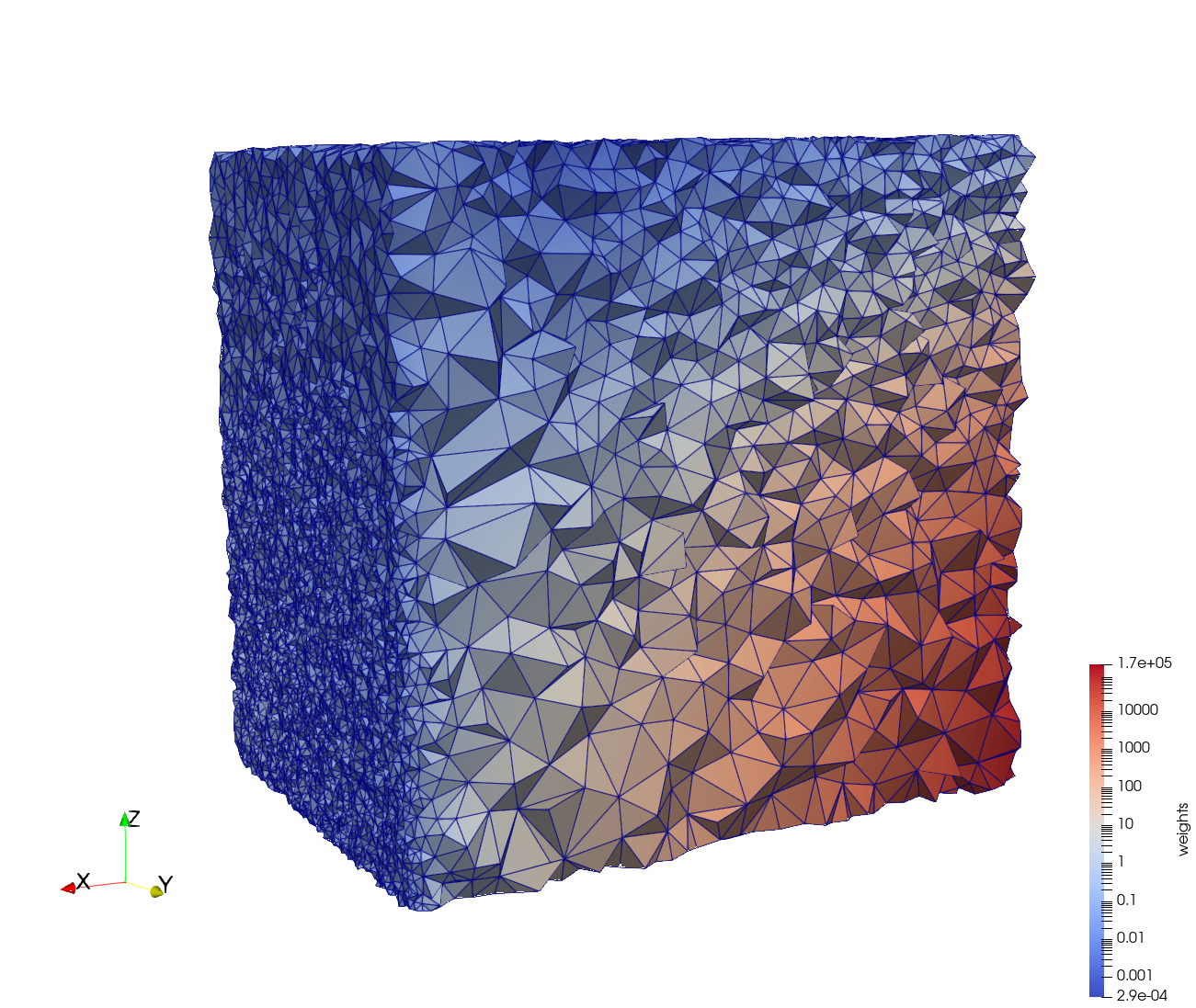}
        \label{fig:adaptative-deviation}
    }
    \caption{(a) The adaptation rule applied on a 2D element, (b) an adaptive tessellation of $Im(H_{u})$ CFF using $SF_1$, (c) an adaptive tessellation of $Im(H_{u})$ CFF using $SF_2$.}
\end{figure}

To efficiently capture the distribution of data, the tessellation needs to \textit{adapt} to the values of the function. The goal of the third approach is to use mesh adaptation to increase the point density of the mesh where it is most needed to reduce the discretization error induced by tessellation. Adaptive tessellations are generated using two different so called \textit{sizing functions}. We use the input image as a background (\textit{BG}) mesh while refining the mesh. Before refinement, the input image is analyzed, and the \textit{image\_weight\_range} $= max(image\_values) - min(image\_values)$ is computed. The refinement algorithm queries the sizing function (\textit{SF}) to verify whether each element satisfies the user constraints. Each time the SF is called upon an element, it creates a set of sampling points \textit{SP} out of the element along with their set of values \textit{V} queried from the \textit{BG}. As Figure \ref{fig:sample-points-quering} indicates, \textit{SP} consists of the element’s vertices, barycenter, and midpoints of vertices and barycenter. The points that are preserved out of \textit{SP} are those that lie within the \textit{BG} mesh. Finally, we compute the \textit{element\_weight\_range} $= max(V) - min(V)$. 

The first version ($SF_1$) of this sizing function  splits an element if $\frac{element\_weight\_range}{image\_weight\_range}$ is greater than a user-defined \textit{weight\_limit}. As Figure \ref{fig:adaptative-gradation} illustrates,  $SF_1$ captures the behavior of the function by creating many simplices where the value of the function is large, which is not necessarily areas where the interpolation error is large.

 The second version ($SF_2$) of this sizing function  splits an element if  $\frac{element\_weight\_range}{|average(V)|}$ is greater than a user-defined \textit{weight\_limit}. As Figure \ref{fig:adaptative-deviation} illustrates, $SF_2$ works such that points are added more densely in regions where the function is rapidly changing, and less points are used in areas where the function is relatively constant.
 
\begin{figure}[!ht]
    \centering
    \subfloat[]{
        \includegraphics[width=0.5\linewidth]{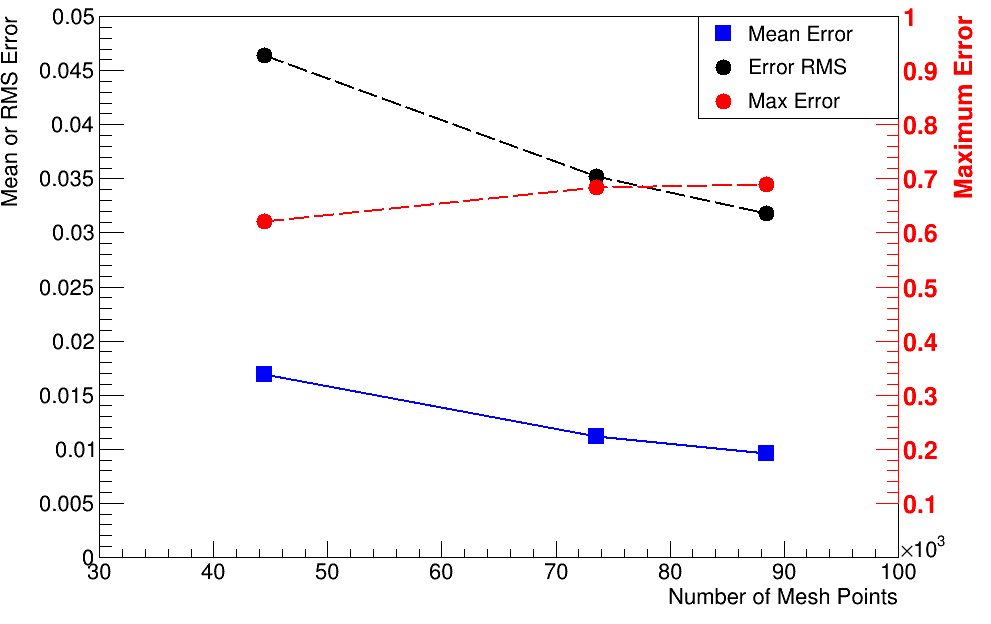}
        \label{fig:error-adaptive-gradation}
    }
    \subfloat[]{
        \includegraphics[width=0.5\linewidth]{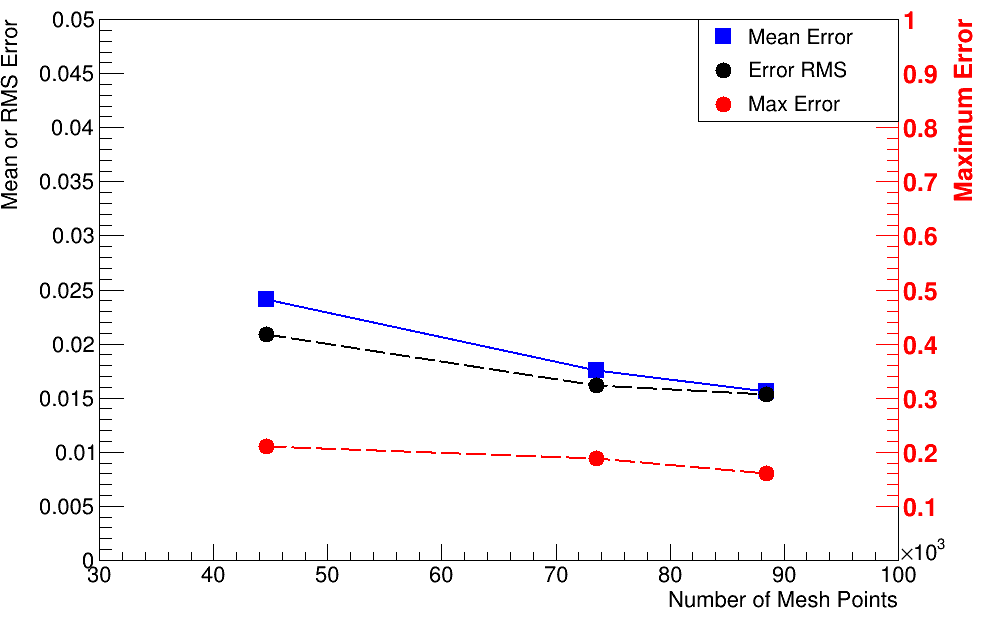}
        \label{fig:error-adaptive-deviation}
    }
    \caption{Mean, RMS, and Max interpolation error statistics of different mesh sizes of the $Im(H_{u})$ CFF using (a) the $SF_1$ and (b) the $SF_2$.}
\end{figure}

 The plot in Figure \ref{fig:error-adaptive-gradation} shows that the $SF_1$ produces a tessellation with a worse error precision compared to a uniform mesh with a similar number of vertices. Figure \ref{fig:error-adaptive-deviation} depicts that the $SF_2$ produces a tessellation that has a lower average interpolation error for a given number of points compared to both the uniform and $SF_1$ methods. However, the maximum interpolation error continues to remain relatively constant.

\subsubsection{Unstructured Iterative-Adaptive Tessellation} \label{sec:iterative-adaptation}

\begin{figure}[H]
    \centering
    \includegraphics[width=0.60\linewidth]{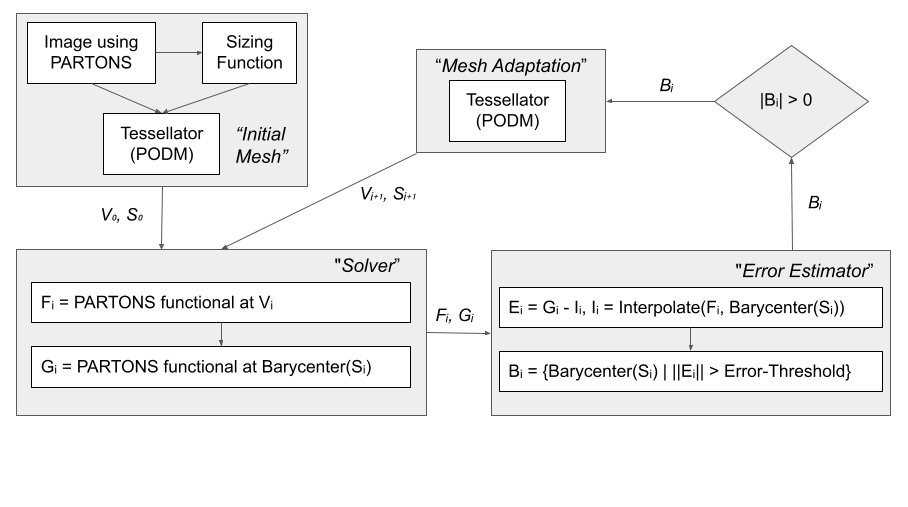}
    \caption{Iterative-adaptive pipeline for CFF data using PARTONS as the solver for error estimation.}
    \label{fig:iterative-adaptivity}
\end{figure}

Computational Fluid dynamics (CFD) applications \cite{roache1997quantification} \cite{fidkowski2011review} have shown that if the mesh is iteratively adapted based on an error estimator, the error can converge much faster than further refining the mesh in a uniform way. Leveraging this idea, Figure \ref{fig:iterative-adaptivity} showcases how a GPD model (\textit{i.e.,} a "solver") can be utilized to define an error estimator for an iterative-adaptive pipeline. In this pipeline, the GPD model (in this case GPDGK16 of PARTONS) is interchangeable and the sizing function is configurable. 

The fourth approach of mesh generation utilizes the interpolation error as a metric for \textit{iterative-adaptation}. As a starting point, an adaptive tessellation is generated using the $SF_2$. Subsequently, for every simplex's barycenter, the value of the function is calculated from both the GPD model and from linear interpolation of the vertices of that simplex. If the interpolation error of a simplex, as defined in section \ref{sec:interpolation-error}, is greater than the user-defined error
threshold, then the simplex's barycenter is added as a new point to the mesh of the next iteration. This iterative process continues until all simplices satisfy the specified error threshold. It should be noted that simplices which satisfy the specified error threshold are not checked again in subsequent iterations. This optimization is of utmost importance since it minimizes the number of calls to the GPD model.

\begin{figure}[!ht]
    \centering
    \subfloat[]{
        \includegraphics[width=0.5\linewidth]{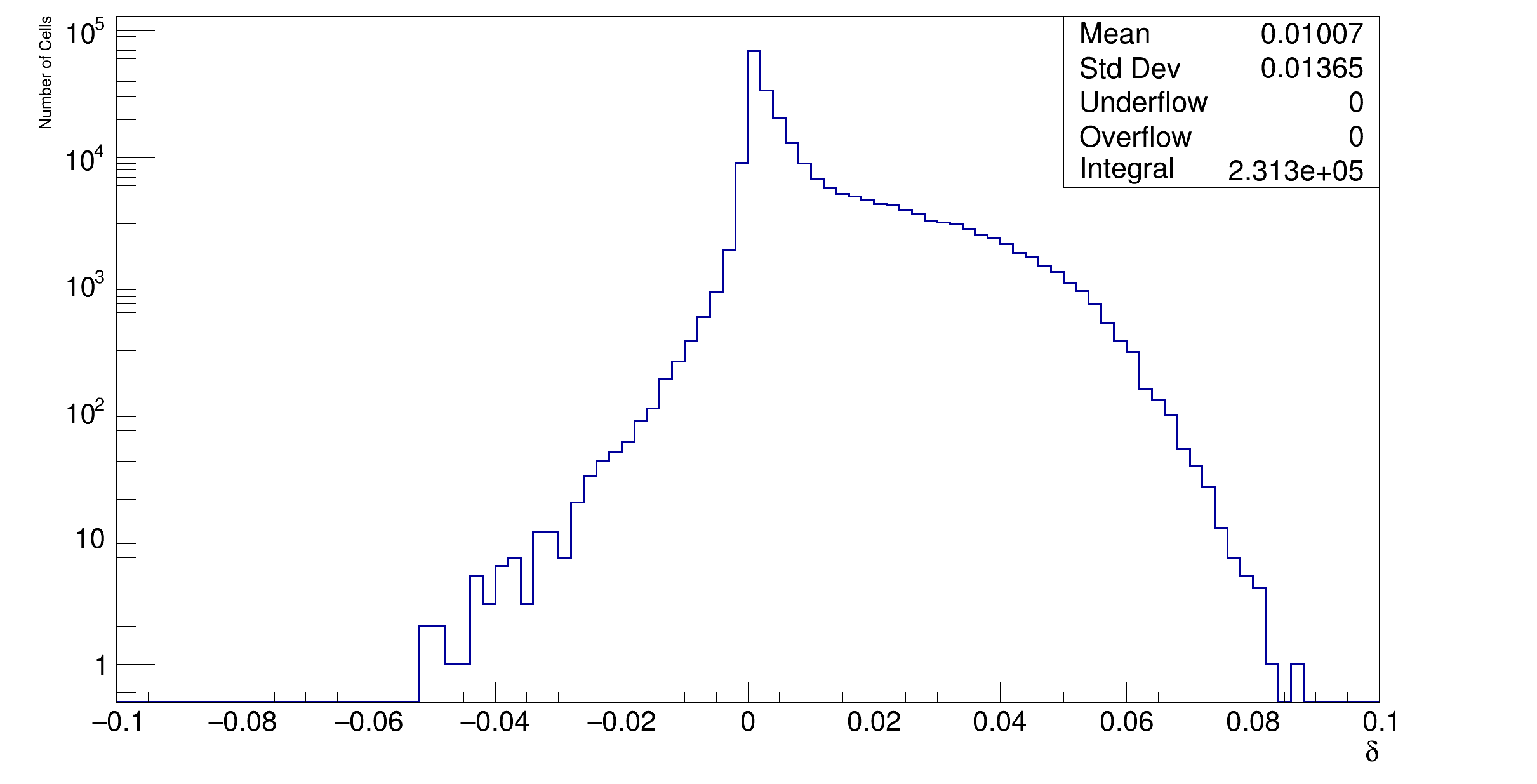}
        \label{fig:ErrorHisto_imag-relative_iteration0}
    }
    \subfloat[]{
        \includegraphics[width=0.5\linewidth]{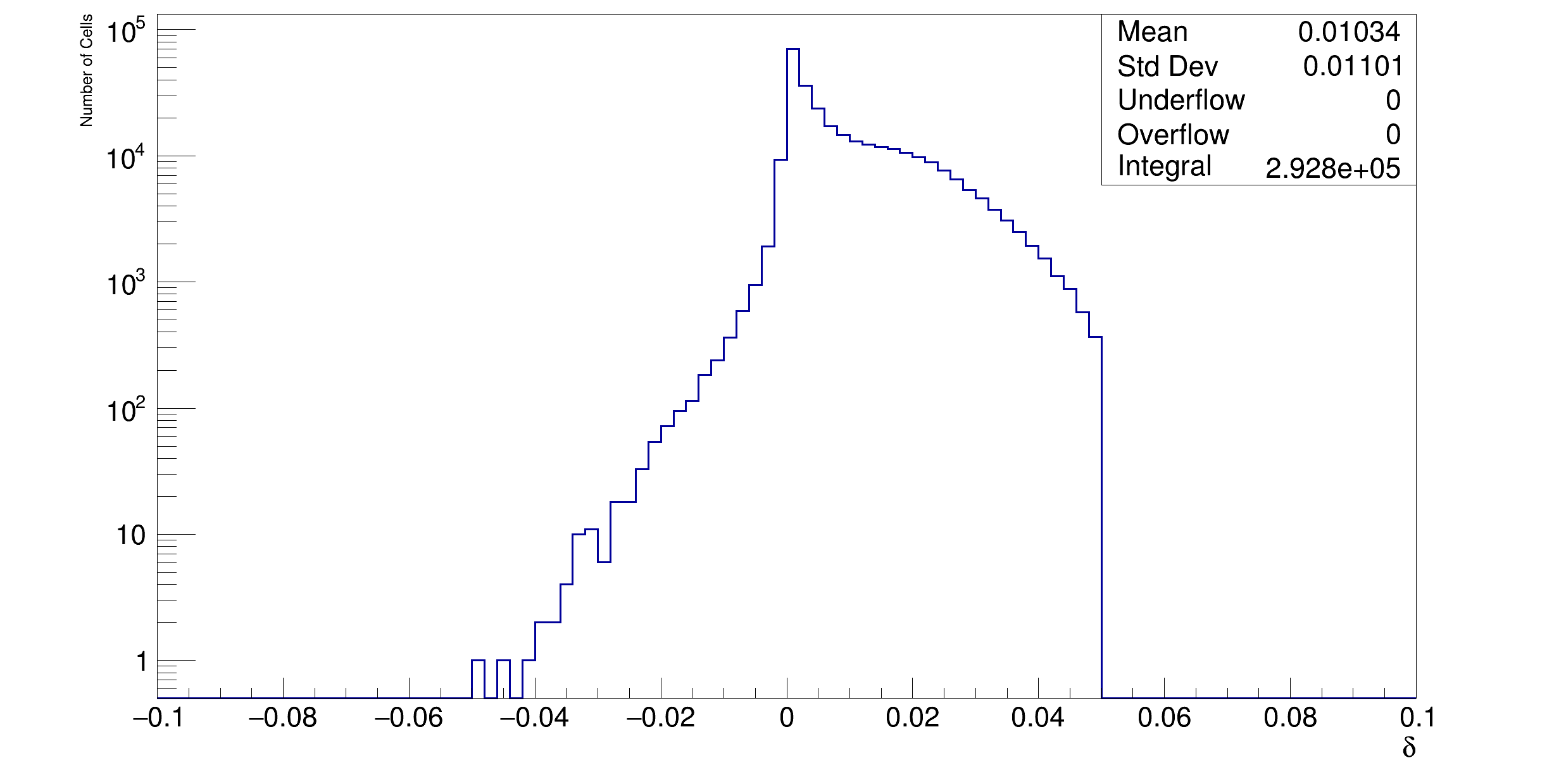}
        \label{fig:ErrorHisto_imag-relative_iteration5}
    }
    \caption{Histograms of the interpolation error of adaptive-iterative meshes of the $Im(H_{u})$ CFF for (a) the 1st iteration and (b) the 6th iteration using the \textit{relative error} metric.}
\end{figure}

Figures \ref{fig:ErrorHisto_imag-relative_iteration0} and \ref{fig:ErrorHisto_imag-relative_iteration5} demonstrate how after $5$ iterations the distribution of interpolation errors not only narrows, but the tails of the distribution also disappear up to the value of the specified error threshold of $0.05$ using the relative error metric. This shows how the iterative method specifically targets simplices with large interpolation error for refinement.

\begin{figure}[!ht]
    \centering
    \subfloat[]{
        \includegraphics[width=0.5\linewidth]{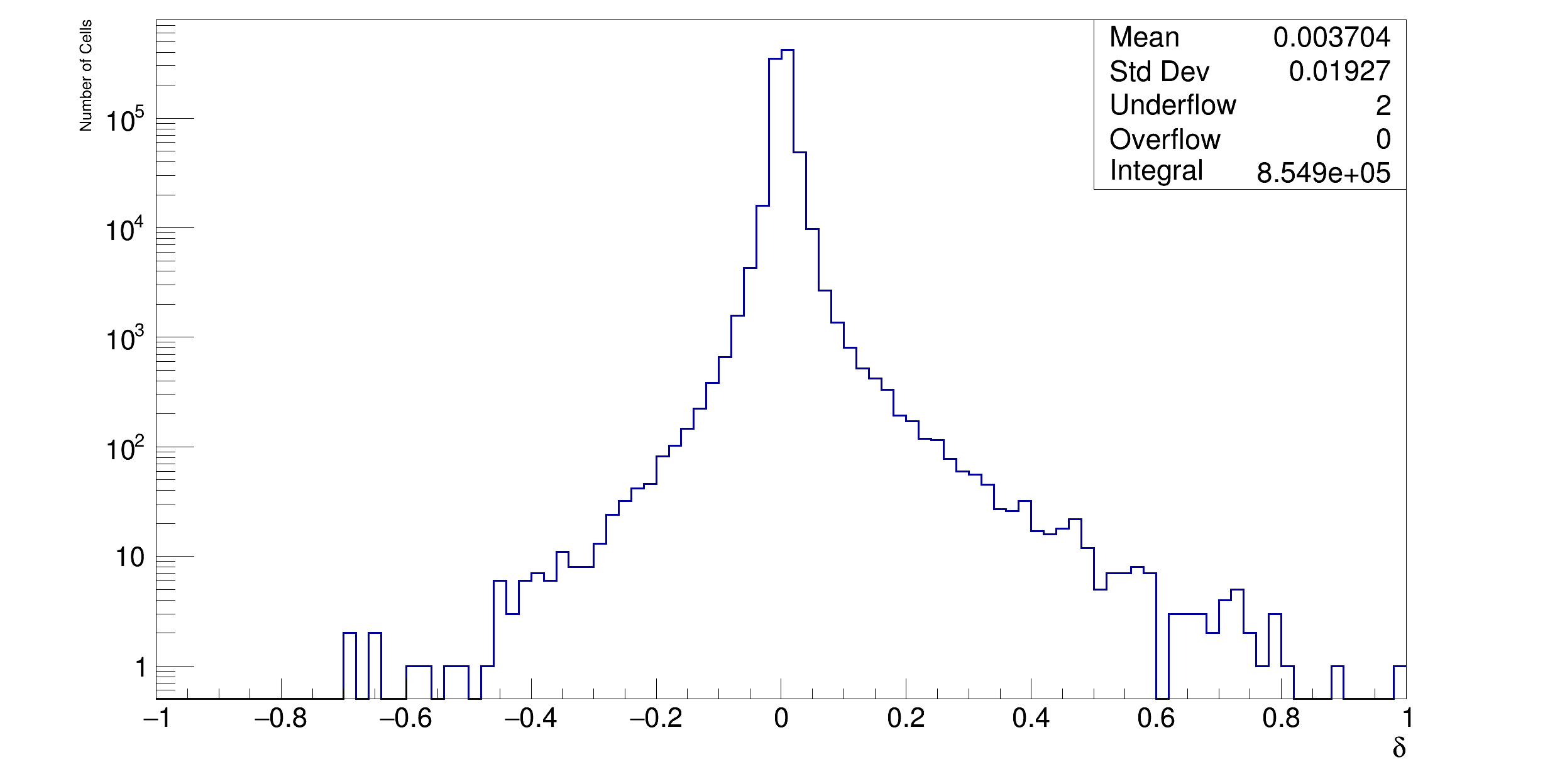}
    }
    \subfloat[]{
        \includegraphics[width=0.5\linewidth]{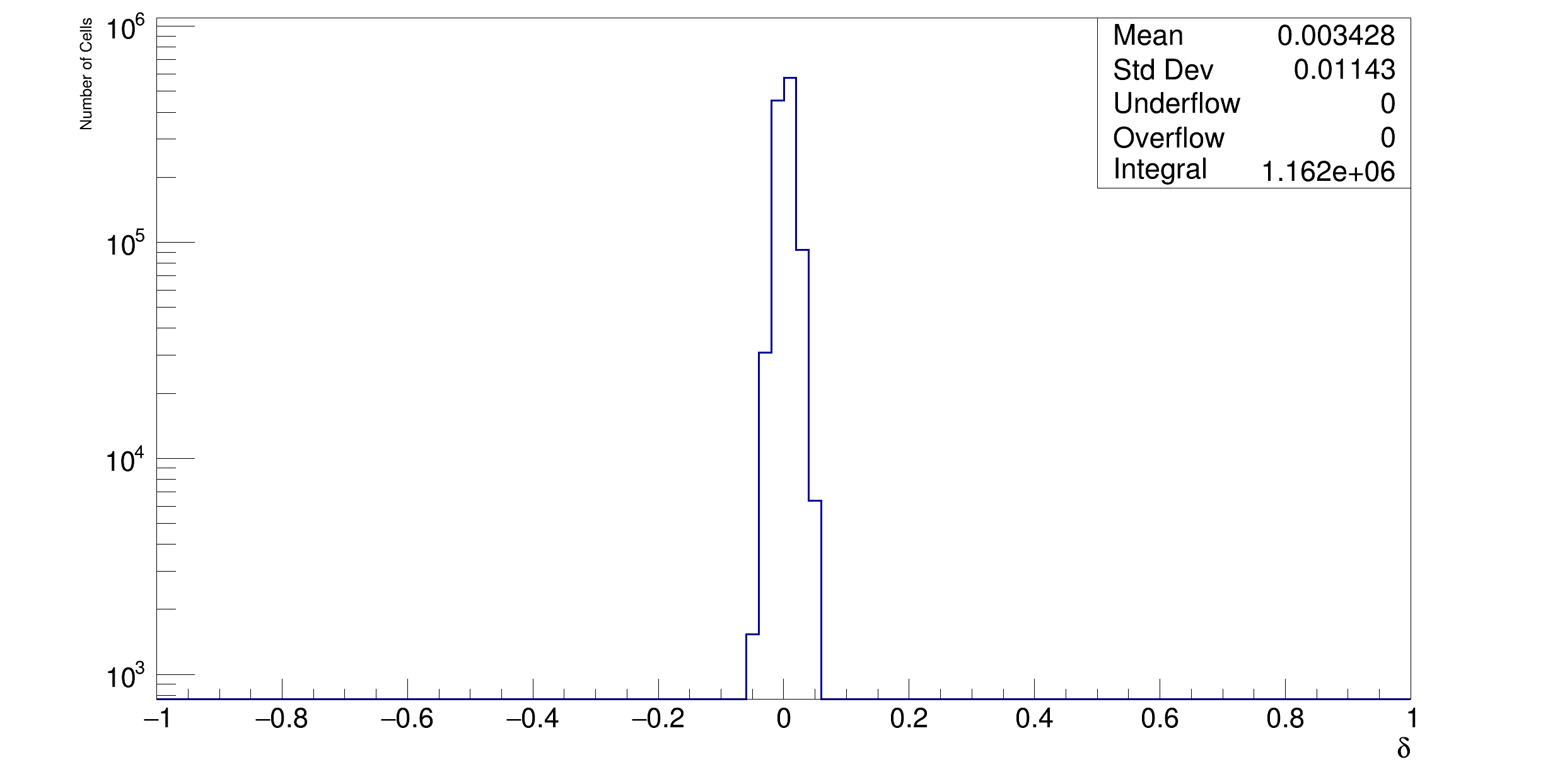}
    }
    \caption{Histograms of the interpolation error of adaptive-iterative meshes of the $Re(H_{u})$ CFF for (a) the 1st iteration and (b) the 11th iteration using the \textit{min(relative, absolute) error} metric.}
    \label{fig:ErrorHisto_real-relative}
\end{figure}

\begin{figure}[H]
    \centering
    \subfloat[]{
        \includegraphics[width=0.5\linewidth]{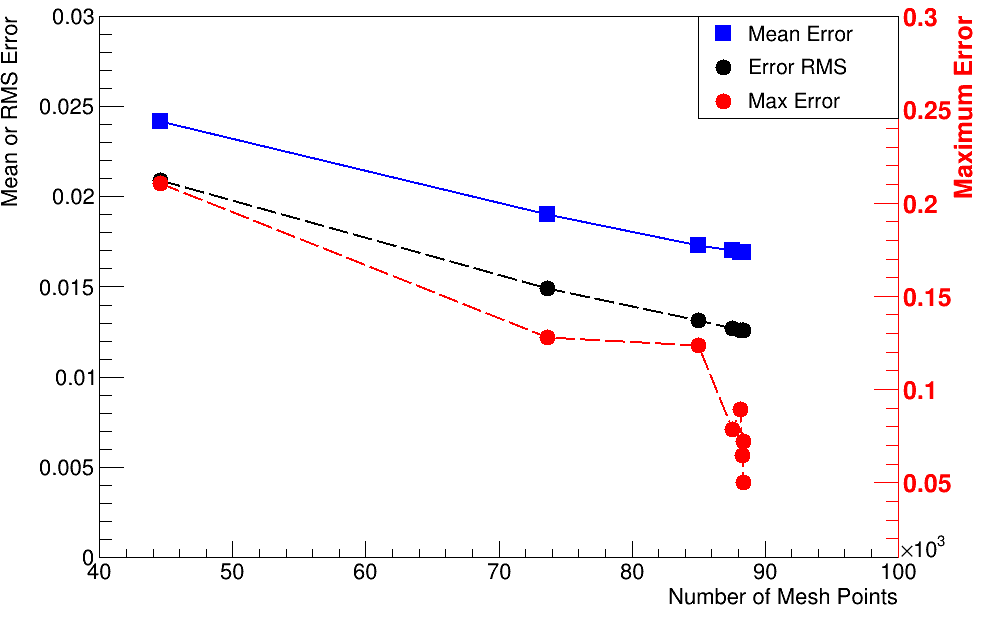}
        \label{fig:meanError_iterative-imag-relative}
    }
    \subfloat[]{
        \includegraphics[width=0.5\linewidth]{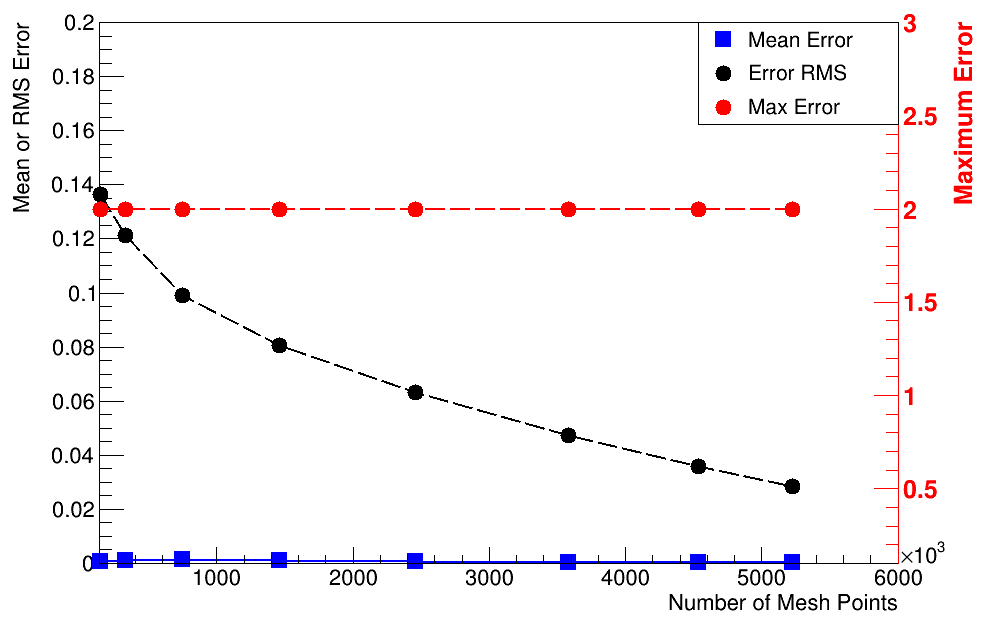}
        \label{fig:meanError_iterative-real-relative}
    }
    \caption{Mean, RMS, and Max interpolation error statistics of different mesh sizes of the (a) $Im(H_{u})$ CFF and (b) $Re(H_{u})$ CFF using the \textit{relative error} metric.}
\end{figure}

\begin{figure}[H]
    \centering
    \subfloat[]{
        \includegraphics[width=0.5\linewidth]{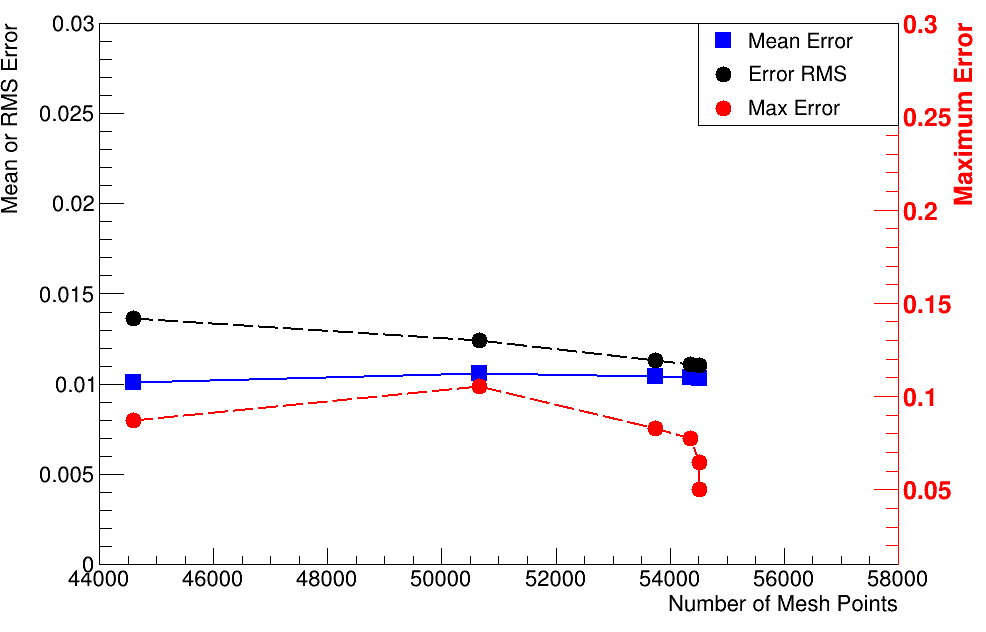}
        \label{fig:meanError_iterative-imag-min}
    }
    \subfloat[]{
        \includegraphics[width=0.5\linewidth]{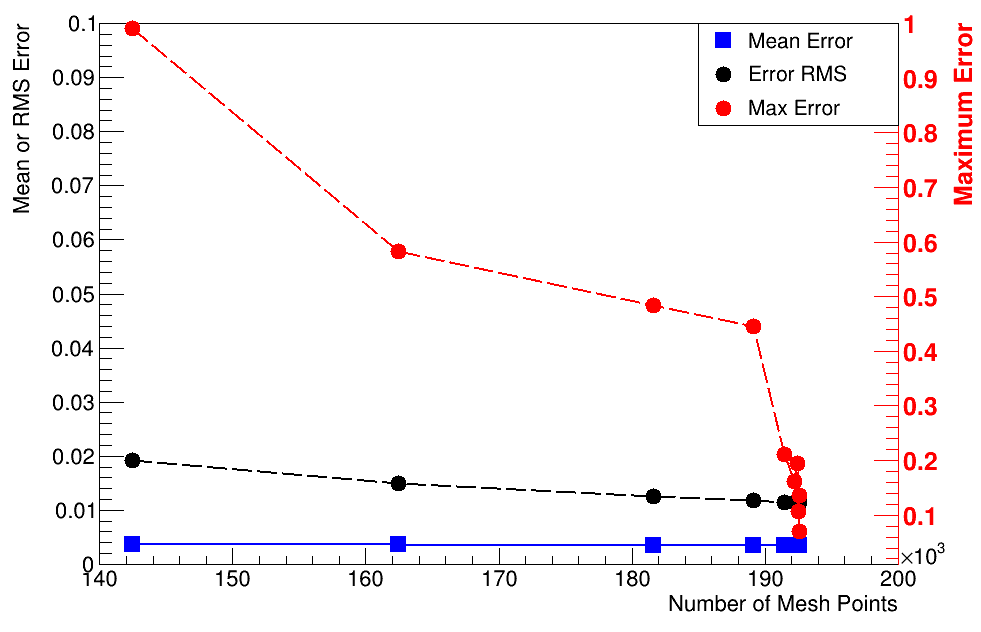}
        \label{fig:meanError_iterative-real-min}
    }
    \caption{Mean, RMS, and Max interpolation error statistics of different mesh sizes of the (a) $Im(H_{u})$ CFF and (b) $Re(H_{u})$ CFF using the \textit{min(relative, absolute) error} metric.}
\end{figure}

Figure \ref{fig:meanError_iterative-imag-relative} shows how this iterative-adaptive process quickly reduces both the mean interpolation error and the maximum interpolation error. However, as explained in section \ref{sec:interpolation-error}, the \textit{relative error} increases significantly when the function has very small positive and negative values. To demonstrate this issue, the iterative process was applied to the $Re(H_{u})$ CFF, which has both positive and negative values, and has very large relative interpolation errors as the function approaches zero. Figure \ref{fig:meanError_iterative-real-relative} shows that this definition of interpolation error will result in a large maximum error, regardless of how many vertices are added (even up to 6 million vertices). This is why it is necessary to define the interpolation error as $min(\textit{absolute error}, \textit{relative error})$. Using this interpolation error as a metric to guide the iterative-adaptation, both functions tested were refined until a maximum interpolation error of 0.05 was met. Figures \ref{fig:meanError_iterative-imag-min} \& \ref{fig:meanError_iterative-real-min} depict how this iterative-adaptation method succeeds in achieving the desired maximum interpolation error. For the $Im(H_{u})$ CFF, the desired maximum interpolation error of 5\% was met in just 6 iterations and 54,500 vertices. Similarly, for the $Re(H_{u})$ CFF, an average error of 0.2\% and maximum error of $5\%$ was met in $11$ iterations and 192,000 vertices. The narrowing of the interpolation error distributions can be seen in Figure~\ref{fig:ErrorHisto_real-relative}.

\subsection{Tessellation utilities} \label{sec:tesselation-utilities}
A tessellation by itself is just a data structure. To access that data, we need certain utilities that are classified into two categories: search and numerical operations.

\subsubsection{Search Operations}
The most common search operation is \textbf{\textit{point location}}. Given a point \textit{p}, and a tessellation with $N_S$ simplices, find the simplex $S_i$, where $i = 0,..., N - 1$, within which \textit{p} is located. This search operation is performed using the \textit{vtkStaticCellLocator} data structure of the VTK library \cite{VTK}. The \textit{vtkStaticCellLocator} includes the functionality to perform the point location operation with worst-case time complexity $O(log(n^3))$, and it is designed to be used in parallel as well.

Another famous search operation is \textbf{\textit{find cell neighbors}}. Given a simplex $S_i$ of a tessellation, find its neighbor simplices with which it shares a face. This search operation is performed using the \textit{vtkStaticCellLinks} data structure of the VTK library \cite{VTK}. This search operation has the time complexity $O(1)$, assuming that every simplex has this information stored inside its data structure, while generating the tessellation.

\subsubsection{Numerical Operations} \label{sec:numerical-operations}

As opposed to \textit{Unstructured} grids, a benefit offered by \textit{Structured} grids is that the computation of the volume of a cell, which is cubic, is trivial. On the other hand, the computation of the \textbf{\textit{volume of a simplex}} is more complex. Let the points $\{\textbf{X}_i|i=0,...,N\}$, which are \textit{column} vectors, define an \textit{N}-simplex \textit{S}. The volume of \textit{S} is expressed with the following formula:

$$
V_S = 
\left|det\begin{pmatrix}
{X_0}^T & 1\\
{X_1}^T & 1\\
... & ... \\
{X_{N}}^T & 1\\
\end{pmatrix}\right| / N!
$$

The volume of a simplex is useful because it can be utilized to perform \textbf{\textit{linear interpolation in a simplex}}. Let the functional values of the points be defined as $\{f_{X_i}|i = 0,...,N\}$, and \textit{p} be an arbitrary point inside \textit{S}. Then, with an arbitrary choice of $X_0$, we define a basis set for \textit{S} and \textit{rhs} column vector as follows:

$$
\{{\Delta X}_i = X_i - X_0 | i = 1, ..., N\},
\\rhs = p - X_0
$$

Afterward, we can define the interpolation weights $\{W_{X_i}|i = 0,...,N\}$ as follows:

$$
\forall X_i | i = 1, ..., N: W_{X_i} = det(M) / V_S, \text{ where } M = N \times N, \text{ and } M_{j} = 
\begin{cases} 
    rhs & j = i \\
    {\Delta X}_i & j \neq i \\
\end{cases}, j = 1, ..., N
$$

% $$
% W_{X_i} = det(M) / V_S, i = 1, ..., N \text{ where } M = N \times N, \text{ and } M_{j} = 
% \begin{cases} 
%     rhs & j = i \\
%     {\Delta X}_i & j \neq i \\
% \end{cases}, j = 1, ..., N
% $$

$$
W_{X_0} = N! - \sum_{i = 1}^{N}W_{X_i}
$$

Utilizing the interpolation weights and the functional values, we can calculate the interpolated functional value of \textit{p} as follows:

$$
f_p = \frac{\sum_{i = 0}^{N} f_{X_i} * W_{X_i}}{\sum_{i = 0}^{N}W_{X_i}}
$$

%%%%%%%%%%%%%%%%%%%%%%%%%%%%%%%%%%%%%%%%%%%%%%%%%
%%%%%%%%%%%%%%%%% Section 4 %%%%%%%%%%%%%%%%%%%%%
%%%%%%%%%%%%%%%%%%%%%%%%%%%%%%%%%%%%%%%%%%%%%%%%%

\section{Monte Carlo Event Generation using Look-Up Table on Tessellated Mesh} \label{sec:monte-carlo-simulation-using-lookup-tables}

Monte Carlo event generators in nuclear femtography often want to generate the final state events according to the probability distribution defined by the differential cross section. To do so, event generators can either generate events uniformly in the phase space and weigh the events by the value of the cross section, or generate unweighted events bound by cross section distribution using an acceptance-rejection method. Both methods require in-situ calculations of the value cross section at every MC point. Additionally, the cross section can vary over several orders of magnitude leading to many rejected events in an acceptance-rejection method.

\begin{figure}[!ht]
    \centering
    \includegraphics[width=0.6\linewidth]{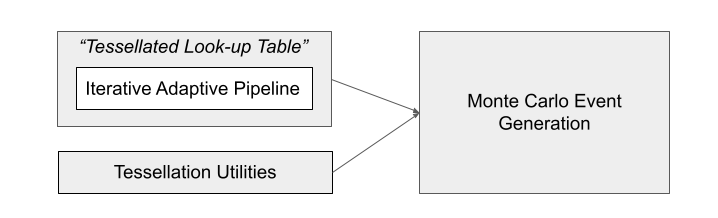}
    \caption{Monte Carlo event generation pipeline using an adaptive tessellated look-up table and the tessellation utilities.}
    \label{fig:Look-up-table-MC}
\end{figure}

An adaptive look-up table of a tessellation of the cross section can be used to eliminate these issues and improve the efficiency of MC event generation.  In the case of DVCS, the 5-fold differential cross section can be reduced to the 4-fold differential unpolarized cross section. Therefore, a 4D tessellation would be needed to be used as a look-up table of the cross section. However, the 3-dimensional Compton Form Factors that contribute to the DVCS cross section can be tessellated and used as an adaptive look-up table in the full calculation of the cross section. Since calculation of the CFFs from models is computationally expensive, querying the tessellation of a CFF as a look-up table will save computational efforts. The proposed Monte Carlo (\textit{MC}) event generation method, which is depicted in Figure \ref{fig:Look-up-table-MC}, can generate events using the adaptive look-up table either uniformly or using an “importance-sampling” algorithm to generate MC events according to the distribution of the CFF. First, the integral of the interpolation function of each simplex is calculated:
 \begin{equation}
        I_{k} = V_{k} \sum_{j=0}^{N} \frac{\bar{f}_{j}(k)}{N+1} 
    \end{equation}
where the volume of each simplex is extracted using the Tessellation utilities.
The Riemann sum of simplex integrals is calculated and normalized to 1. 
\begin{equation}
        F_{0}=0, \qquad F_{m} = \frac{\sum_{k=1}^{m} I_{k}}{\sum_{k=1}^{N_{s}} I_{k}} 
\end{equation}
Then a random MC event $u\in[0,1]$ is generated. The simplex that is selected, using binary search, has the random event weight bounded by its Riemann sum, such that:
\begin{equation}
    F_{k-1} \leq u < F_{k}
\end{equation}

The MC event, called $P_{nS}$, is randomly generated using a uniform distribution \cite{uniform-picking} in the selected \textit{n}-dimensional simplex, which has $X_i$ as vertices. The weights that are used are denoted as $\lambda_j, j=0,...,n+1$.

\begin{equation}
\begin{split}
    P_{nS} = \sum_{i=1}^{n+1} \left( (1-\lambda_i) \prod_{j=0}^{i-1}\lambda_j\right) X_i, \text{ with }\lambda_0 = 1\text{ and } \lambda_{n+1} = 1\\
    \forall\lambda_j | i = 1, ..., n: \lambda_j = \sqrt[k]{z_j}, \text{ with } k = (n+1) - j
\end{split}
\end{equation}

\begin{figure}[!ht]
    \centering
    \subfloat[]{
        \includegraphics[width=0.45\linewidth]{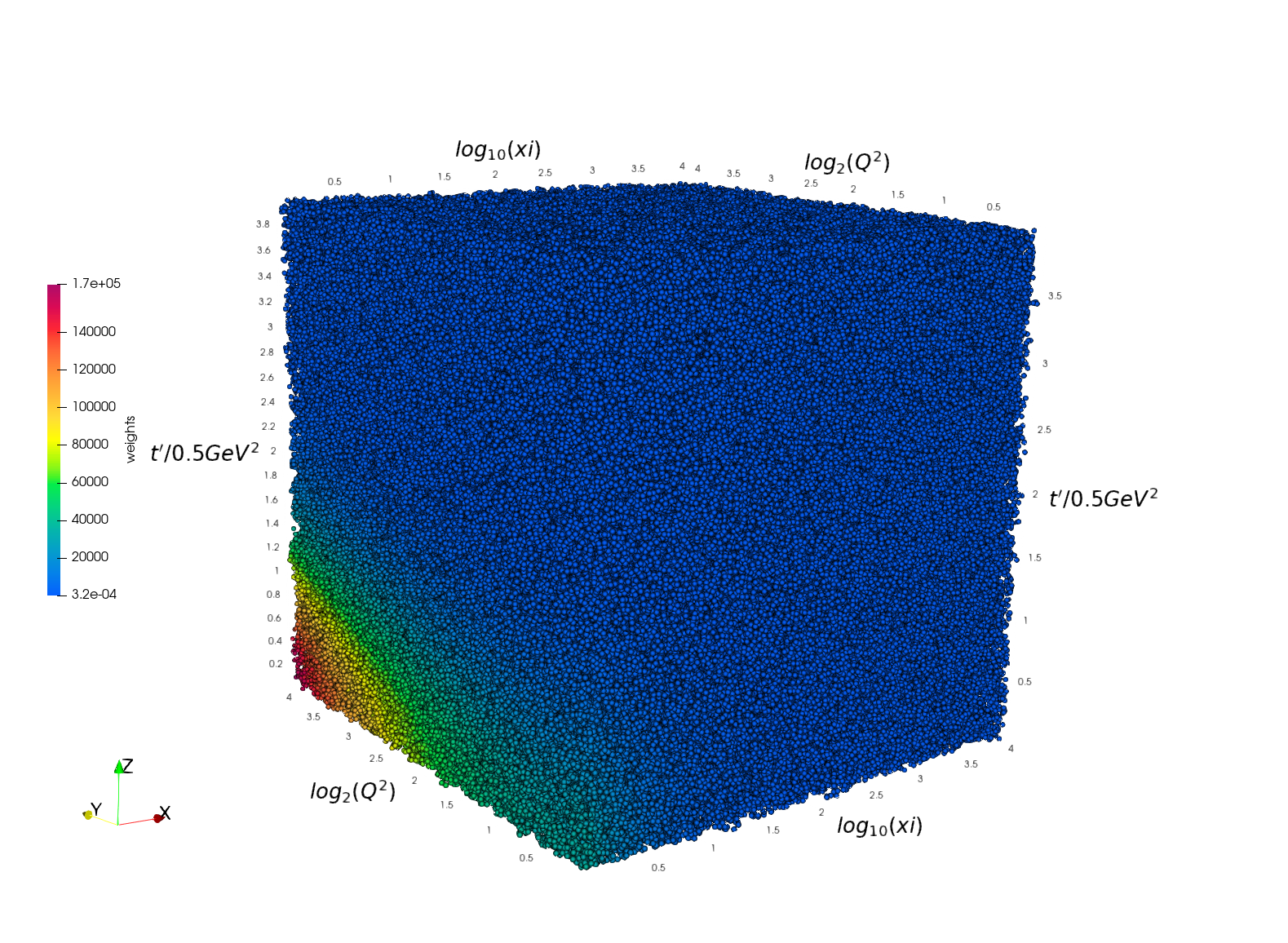}
    }
    \subfloat[]{
        \includegraphics[width=0.45\linewidth]{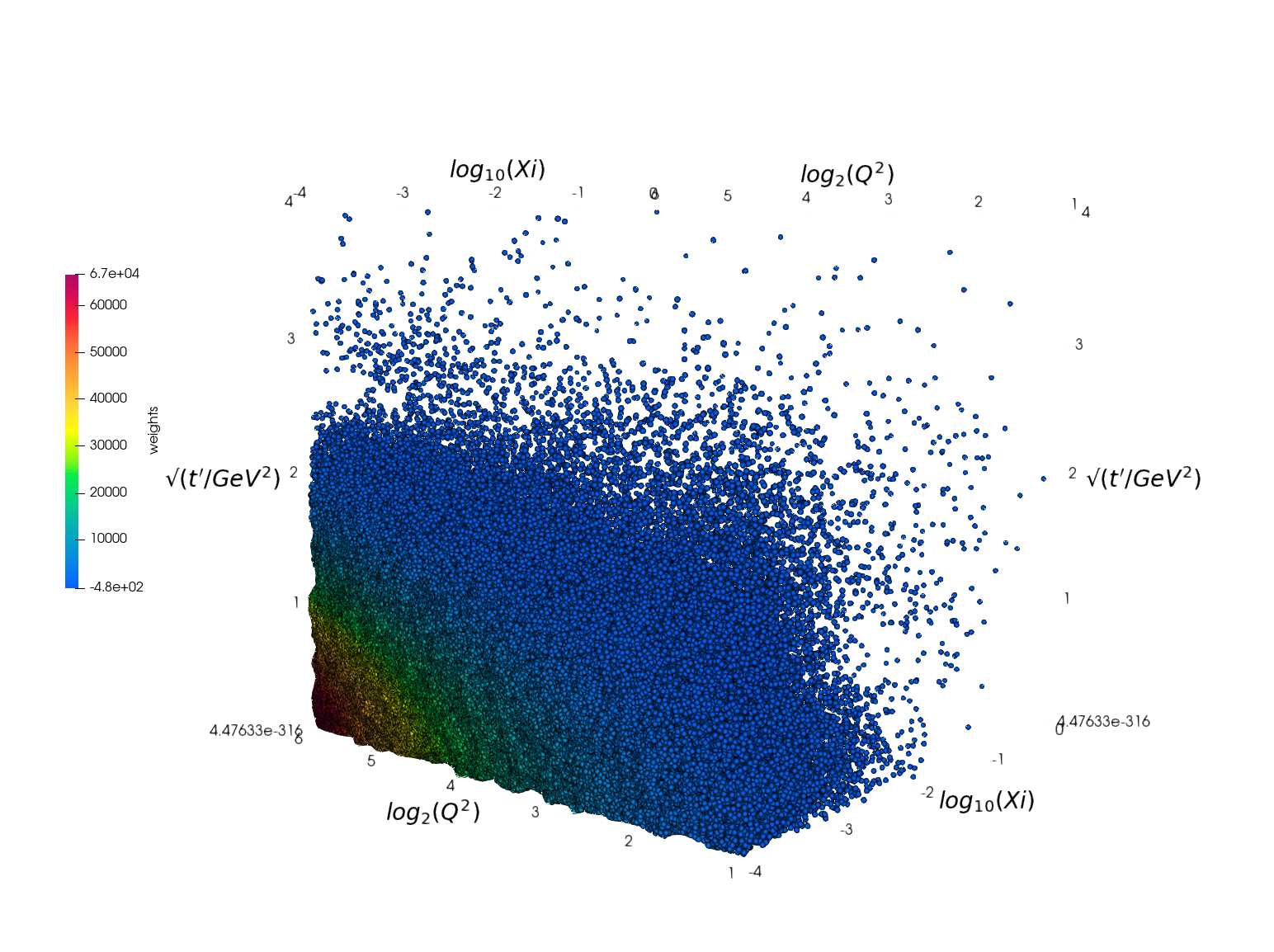}
    }
    \caption{1 million Monte Carlo events generated (a) uniformly and (b) according to the distribution of the $Im(H_{u})$ CFF using an adaptive tessellation as a lookup table.}
    \label{fig:event-generation}
\end{figure}

When performing full detector simulations for an experiment, the Monte Carlo events are often weighted by the cross section or generated according to the probability distribution defined by the cross section. However, this case study created 3-D tessellations of CFFs rather than the 5-D cross section. As a proof of concept, the tessellation utilities were used to extract the CFFs from a tessellation and generate Monte Carlo events according to the probability distribution defined by the CFF.
This method of Monte Carlo event generation using an adaptive look-up table to extract the CFF from a tessellation significantly increases efficiency. Figure \ref{fig:event-generation} shows a visualization of 1 million MC events generated according to the distribution of the $Im(H_{u})$ Compton Form Factor. As mentioned previously, the drawback of this method is that the MC events will have a small error from the linear interpolation. The significance of the computation time saved and the performance optimizations of this method are explained in section \ref{sec:monte-carlo-performance}.

%%%%%%%%%%%%%%%%%%%%%%%%%%%%%%%%%%%%%%%%%%%%%%%%%
%%%%%%%%%%%%%%%%% Section 5 %%%%%%%%%%%%%%%%%%%%%
%%%%%%%%%%%%%%%%%%%%%%%%%%%%%%%%%%%%%%%%%%%%%%%%%

\section{Performance Evaluation}
This performance evaluation consists of four sections. The first section shows the total run-time of the iterative-adaptive pipeline to generate a mesh of $Im(H_{u})$ using (a) 40 cores, and (b) 800 cores for PARTONS, and the run-time per iteration of the iterative-adaptive pipeline steps to generate a mesh of (a) $Im(H_{u})$, and (b) $Re(H_{u})$ using 800 cores for PARTONS. The second section estimates the cost of generating an unstructured uniform tessellation of $Im(H_{u})$ that satisfies the error threshold, and compares it with the resulting tessellation of the iterative-adaptive pipeline using 800 cores for the PARTONS. The third section showcases the benefits of generating MC events utilizing the tessellation of $Im(H_u)$ as a look-up table generated by the iterative-adaptive pipeline, compared to the GPD model using (a) 40 cores, and (b) 800 cores. All iterative-adaptive tessellations are generated using the \textit{min(relative, absolute) error} metric and error threshold = 0.05.

To evaluate the performance of our results, we use the CPU nodes of Old Dominion University's Wahab HPC cluster that each have two 20-core Intel Xeon Gold 6148 2.40GHz CPUs (totaling 40 cores per node) and 384GB of RAM. The compiler that we use is GCC 7.5.0.

\subsection{Run-time of Iterative-Adaptive Pipeline steps}

\begin{table}[H]
    \centering
    \caption{Iterative-Adaptive Pipeline steps.}
    \begin{tabular}{|p{0.09\linewidth}|p{0.75\linewidth}|}
        \hline
        \textbf{Pipeline Step} & \textbf{Description} \\ \hline
        Meshing IO & The Meshing IO step includes the reading and writing of files needed for the Meshing step. \\ \hline
        Meshing & The Meshing step includes the generation of the first mesh and the insertion of additional points to perform adaptive mesh refinement. \\ \hline
        Kinematics & The Kinematics step includes the generation of input, such as initial grid points, mesh vertices, and barycenters, for the PARTONS step. \\ \hline
        PARTONS & The PARTONS step includes the computation of all the CFF values of the input kinematics points, such as input image grid points, mesh vertices, and the mesh simplices' barycenters. \\ \hline
        PARTONS IO & The PARTONS IO step includes the merging of the PARTONS result of each processing core into one file. \\ \hline
        Update Weights & The Update Weights step includes the updating of the weights (CFF values) of the new mesh vertices.\\ \hline
        Compute Error & The Compute Error step includes the linear interpolation computation (using the tessellation utilities) of all the barycenters, the calculation of the linear interpolation error, and the writing of the points that will be added into the next iteration's mesh.\\ \hline
    \end{tabular}%
    \label{tab:iterative-adaptive-pipeline-steps}
\end{table}

\begin{figure}[H]
    \centering
    \subfloat[]{
        \includegraphics[width=0.5\linewidth]{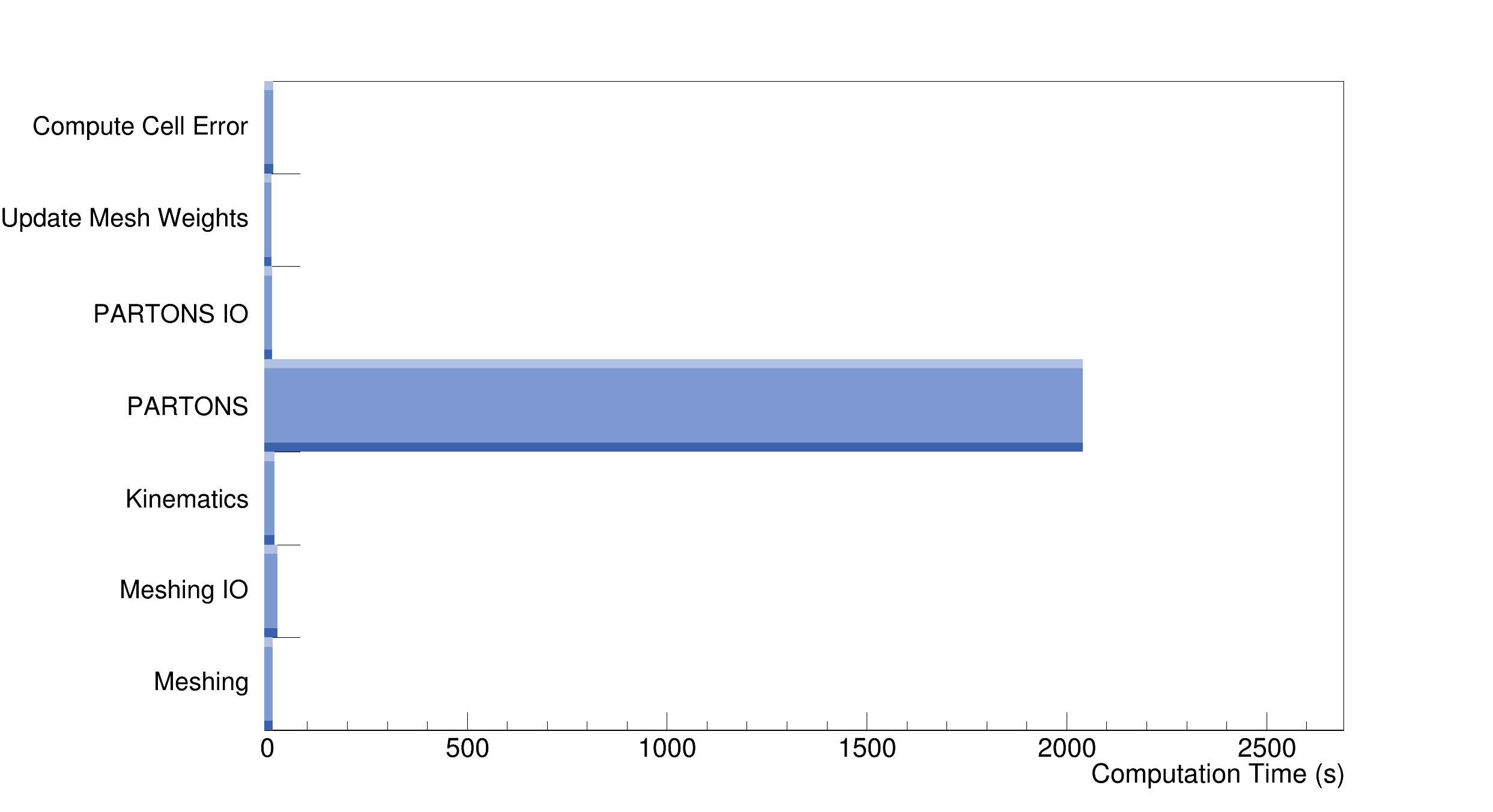}
    }
    \subfloat[]{
        \includegraphics[width=0.5\linewidth]{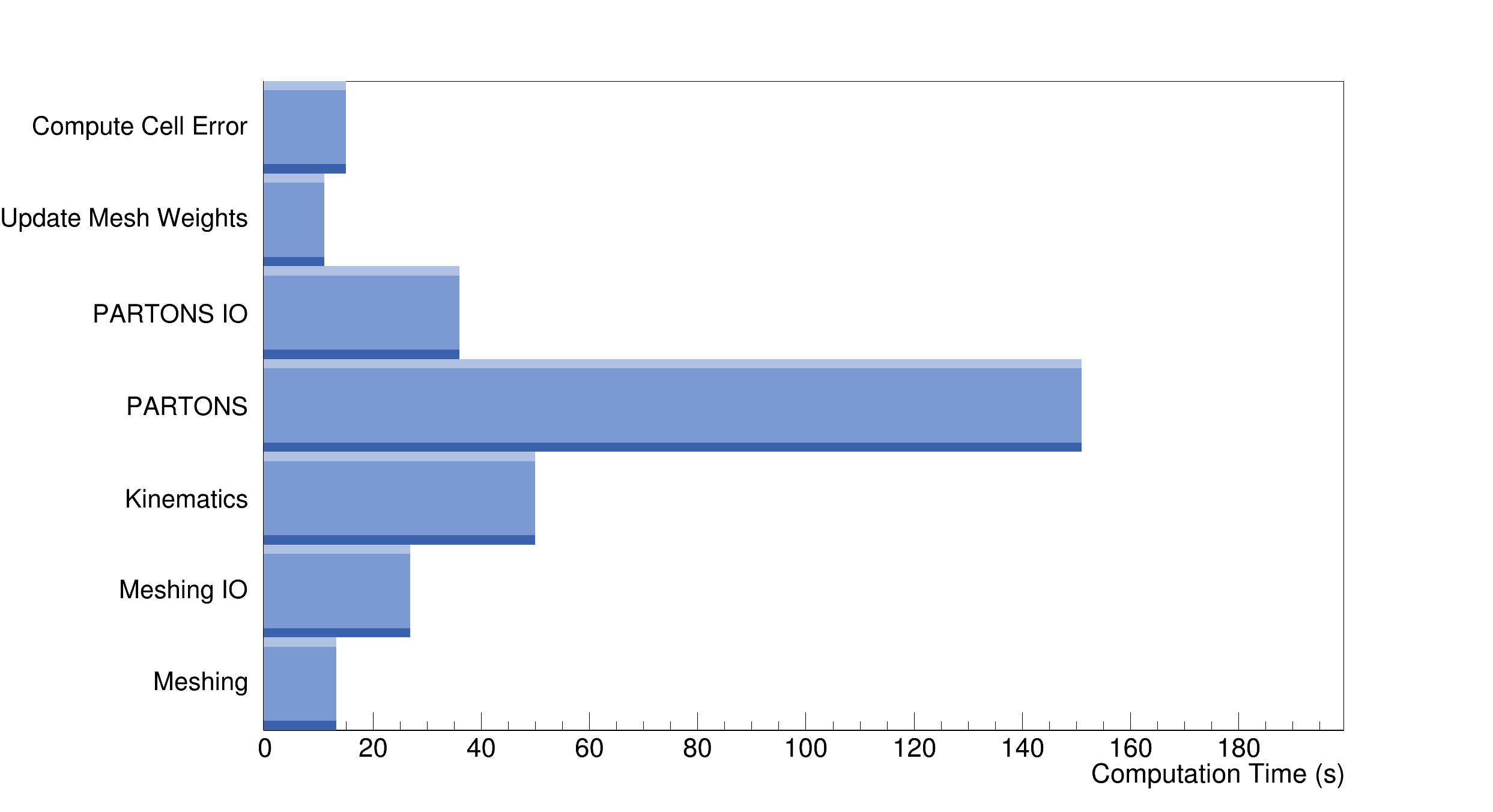}
    }
    \caption{The total run-time of the steps in the iterative-adaptive pipeline to generate a tessellation of $Im(H_{u})$ using (a) 40 cores, and (b) 800 cores for PARTONS.}
    \label{fig:pipeline-steps-total-time}
\end{figure}

The iterative-adaptive pipeline that generates an accurate tessellation for a specific CFF consists of seven core steps as described in Table ~\ref{tab:iterative-adaptive-pipeline-steps}. It should be noted that the step of converting the PARTONS data to the initial image is omitted because it does not constitute a \textit{full} iteration, and its cost (time-wise) is negligible. The meshing step, which utilizes the parallel software PODM, is not executed in parallel. The reasoning behind this choice is that PODM's operation of inserting additional points, which has been implemented specifically for this application, is not yet parallel. The PARTONS step has been modified to effectively utilize distributed computing nodes for parallel execution. As Figure \ref{fig:pipeline-steps-total-time} illustrates, parallelizing the PARTONS step was the most significant factor in reducing the total run-time of the pipeline. The remainder of the pipeline steps are sequential.

\begin{figure}[H]
    \centering
    \subfloat[]{
        \includegraphics[width=0.5\linewidth]{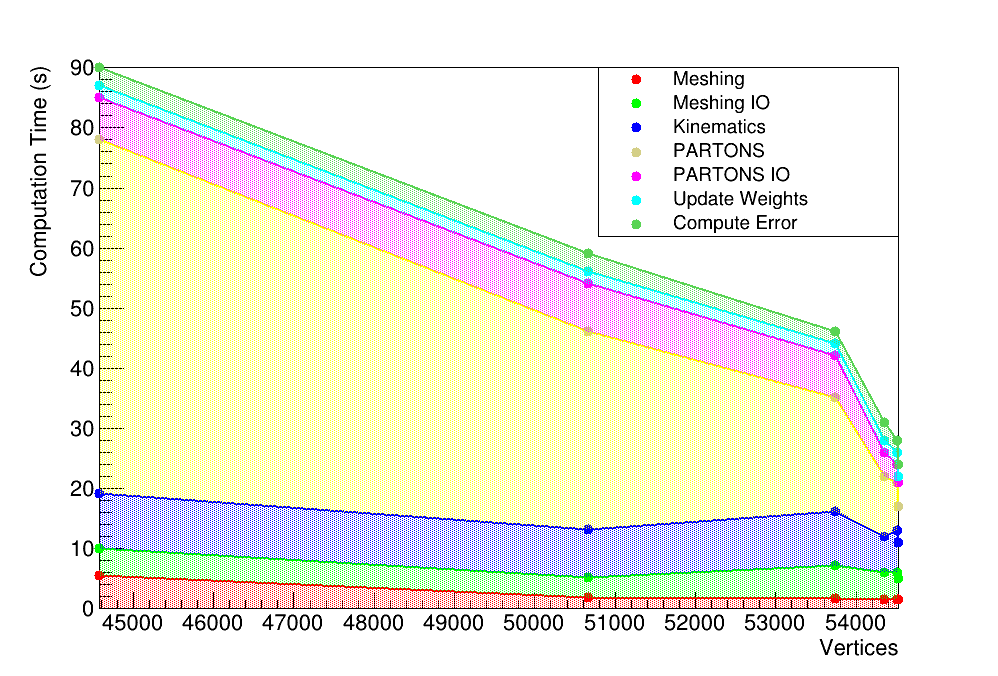}
    }
    \subfloat[]{
        \includegraphics[width=0.5\linewidth]{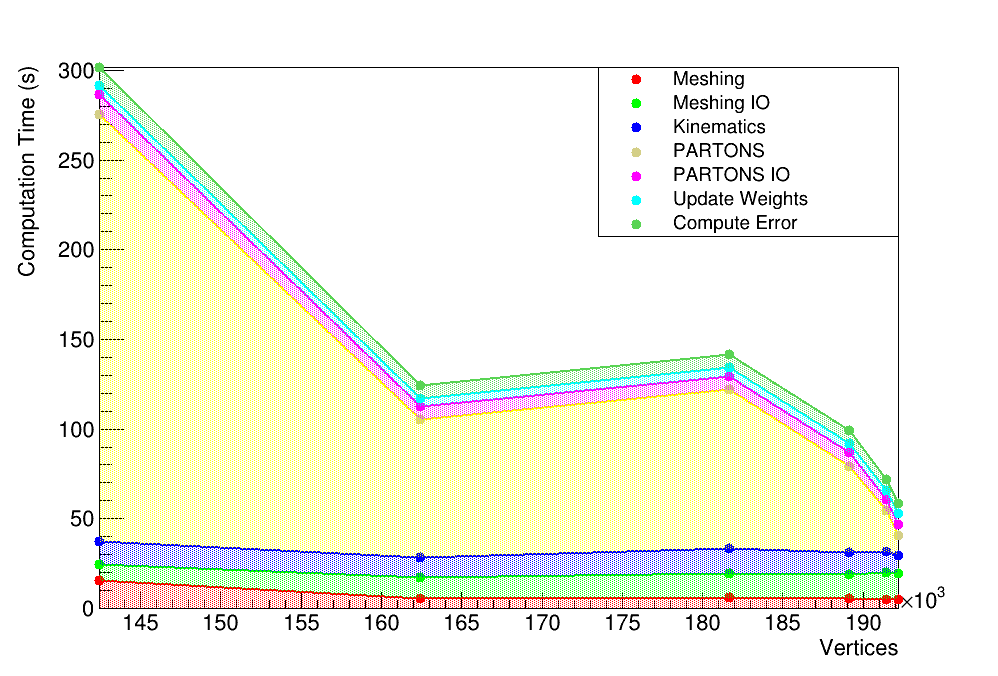}
    }
    \caption{Shown is the run-time of each step per iteration in the iterative-adaptive pipeline to generate a mesh of (a) $Im(H_{u})$, and (b) $Re(H_{u})$ with a maximum interpolation error of 0.05. The widths of the individual bands correspond to the computational time of each step. The PARTONS step was executed using 800 cores.}
    \label{fig:800core-steps-per-iteration-time}
\end{figure}

In section \ref{sec:iterative-adaptation}, we discussed the benefit of ensuring that the iterative-adaptive pipeline knows when a simplex has already satisfied the specified error threshold in a previous iteration. Consequently, in each iteration, only the newly introduced simplices should be checked to see if they satisfy the specified error threshold. Figure \ref{fig:800core-steps-per-iteration-time} verifies our argument, since it showcases that while the number of points (and therefore simplices) increases per iteration, the time spent for PARTONS calls does not follow the same trend and decreases in almost every iteration. 

% \begin{figure}[H]
%     \centering
%     \subfloat[]{
%         \includegraphics[width=0.5\linewidth]{./Figures/Section-5/800imag-min-cumulative.png}
%     }
%     \subfloat[]{
%         \includegraphics[width=0.5\linewidth]{./Figures/Section-5/800real-min-cumulative.png}
%     }
%     \caption{Cumulative run time at each step in the iterative-adaptive pipeline to generate a mesh of (a) $Im(H_{u})$, and (b) $Re(H_{u})$ with a maximum interpolation error of 0.05. The widths of the individual bands correspond to the computational time of each step. The "PARTONS" step was ran on 800 cores.}
%     \label{fig:800core-steps-per-iteration-time-cumulative}
% \end{figure}

\subsection{Uniform Unstructured vs Iterative-Adaptive Unstructured}

In section \ref{sec:unstructured-uniform}, we understood that generating a uniform unstructured tessellation, which satisfies the specified error threshold, is a task with several disadvantages. In this section, we attempt to estimate, using extrapolation techniques, the number of vertices and the computational cost that would be required to generate an accurate uniform unstructured tessellation of $Im(H_{u})$ that satisfies the error threshold using 800 cores for the PARTONS. Subsequently, we compare this computational cost with that of the resulting tessellation of the iterative-adaptive pipeline using 800 cores for the PARTONS.

\begin{table}[H]
    \centering
    \caption{Maximum interpolation error and computation cost of generating uniform unstructured tessellations of $Im(H_{u})$ with different numbers of vertices.}
    \label{tab:estimation-uniform-data}
    \resizebox{0.45\linewidth}{!}{
    \begin{tabular}{|c|c|c|}
        \hline
        \textbf{Vertices} & \textbf{Max Error} & \textbf{Computation cost (sec)} \\ \hline
        44585 & 0.2066 & 204 \\ \hline
        73616 & 0.1451 & 319 \\ \hline
        88388 & 0.1129 & 367 \\ \hline
    \end{tabular}%
    }
\end{table}

Based on the data provided in Table~\ref{tab:estimation-uniform-data}, we derived the following approximation functions, using linear regression, that capture the trend of the maximum interpolation error and the computation cost: $Error_{max}(v) = -2.14*10^{-6}*v + 0.302$, and $ComputationCost(v) = 3.75*10^{-3}*v + 38.1$. Since the desired error threshold is 0.05, the estimated required number of vertices is 117,757, and the estimated required computation cost is 479 seconds. The iterative-adaptive pipeline generated a tessellation with 54,523 vertices in 309 seconds; therefore, it is better in all regards compared to the estimated uniform tessellated result.

\subsection{Monte Carlo Event Generation} \label{sec:monte-carlo-performance}

As described in section \ref{sec:monte-carlo-simulation-using-lookup-tables}, tessellated meshes can be utilized as look-up tables for Monte Carlo (MC) event generation. Consequently, we developed a parallel software for MC simulation capable of utilizing distributed computing nodes. The purpose of this software is to demonstrate that the MC simulation can be greatly expedited using an accurate adaptive tessellation of a CFF. To validate our assumption, we (a) generate an accurate tessellation of $Im(H_{u})$ using the iterative-adaptive pipeline, (b) generate $10^7$ uniformly distributed MC events using the tessellation of (a), and (c) calculate the CFF values of the 10 million MC events generated in (b), using 40 and 800 cores.

\begin{figure}[H]
    \centering
    \includegraphics[width=0.5\linewidth]{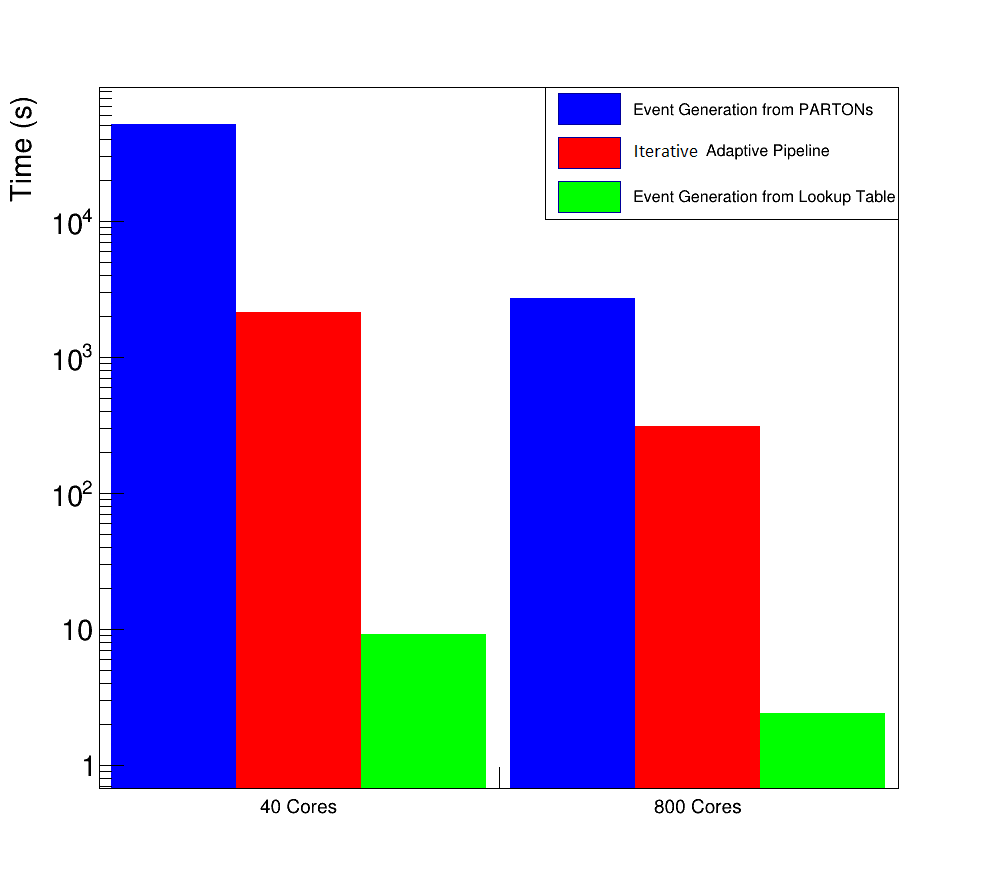}
    \caption{Time comparison to complete the adaptive mesh generation pipeline (red), the time to generate 10 million events by using a mesh as a smart lookup table (green), or directly from PARTONS (blue).}
    \label{fig:MC-timing}
\end{figure}

As Figure \ref{fig:MC-timing} illustrates, the cost to generate an accurate tessellation of $Im(H_{u})$ and $10^7$ events using our MC simulation software is 2,151 and 318 seconds using 40 and 800 cores, respectively. The cost to calculate the CFF values of these $10^7$ events using just PARTONS is 50,997 and 2694 seconds using 40 and 800 cores, respectively. Consequently, the speedup is 23.7 for 40 cores, and 8.47 for 800 cores. In summary, using 40 cores would take about 35.9 minutes instead of 14.2 hours, and using 800 cores would take about 5.3 minutes instead of 44.9 minutes.

A key take-away is that if someone wanted to generate 10 million MC events again, the cost would be dramatically reduced because the tessellation would have already been generated. Based on the aforementioned performance data, we can derive an estimated generalization of the speed-up formulas for 40 and 800 cores, with respect to the number of events $n_e$, as shown in Equations \ref{eq:total-speedup-general} and \ref{eq:total-speedup-specific}.

\begin{equation}
S_{c}(n_e) = \frac{n_e*T_{{PARTONS}_{c}}(1)}{T_{{Iterative\_adaptive\_pipeline}_{c}} + n_e * T_{MC_{c}}(1)} \approx \frac{n_e*\frac{T_{{PARTONS}_{c}}(10^7)}{10^7}}{T_{{Iterative\_adaptive\_pipeline}_{c}} + n_e * \frac{T_{MC_{c}}(10^7)}{10^7}}
\label{eq:total-speedup-general}
\end{equation}

\begin{equation}
S_{40}(n_e) \approx \frac{n_e*\frac{50,997}{10^7}}{2,142 + n_e * \frac{9.2}{10^7}}, \hspace{1cm} S_{800}(n_e) \approx \frac{n_e*\frac{2,694}{10^7}}{309 + n_e * \frac{2.4}{10^7}}.
\label{eq:total-speedup-specific}
\end{equation}

Using these formulas, we estimate that the speedup for generating $10^{10}$ events would be 4,496 and 994 for 40 and 800 cores, respectively. In summary, one might spend about 3.2 hours instead of 19 months when utilizing 40 cores, or 45.2 minutes instead of 1 month when utilizing 800 cores.

\section{Conclusions}
In this work, we designed an iterative-adaptive pipeline capable of efficiently generating an adaptive tessellation for a specific Compton Form Factor of a GPD model with a minimal number of vertices while still maintaining a 1\% mean interpolation error and a 5\% maximum interpolation error. Using such a tessellation as a lookup table significantly decreases the computation time for Monte Carlo event generation, by about 23x for $10^{7}$ events (and using extrapolation, by about 955x for $10^{10}$ events). Equally as important is the ability to re-use a tessellation multiple times for different analyses without having to do any calculations from a physics model.

In the future, since a GPD model has multiple CFFs, we plan to extend this iterative-adaptive pipeline to generate a \textit{common denominator} adaptive tessellation which satisfies the error threshold for all the CFFs. Generating one tessellation (instead of many) will be beneficial for two reasons: 1) it minimizes the GPD calculations, since PARTONS (the GPD model that we used) calculates all the CFFs at once and 2) MC event generation using multiple tessellations as a lookup table would require the usage of the point location operation (which has $O(log(n^3))$ time complexity) to detect a randomly generated point in another CFF tessellation. The common denominator adaptive tessellation of all the 3-dimensional CFFs, as mentioned in section \ref{sec:monte-carlo-simulation-using-lookup-tables}, will be the foundation for efficiently calculating the full DVCS cross section. A potential alternative is to use this iterative-adaptive pipeline to generate a 4-dimensional tessellation of the unpolarized cross section directly from a GPD model. For this purpose, the 4D version of PODM \cite{PODM4D} could be enhanced to work with the iterative-adaptive pipeline. 

These software tools will be made publicly available to the nuclear physics community in hopes that they are utilized to make a database of tessellation lookup tables of Compton Form Factors and other functions from various physics models.

\section{Acknowledgments}
This material is based upon work supported by the U.S. Department of Energy, Office of Science, Office of Nuclear Physics under contract DE-AC05-06OR23177. This work was also funded in part by NSF grant No. CCF-1439079, the Richard T. Cheng Endowment, and the National Institute of General Medical Sciences of the National Institutes of Health under Award Number 1T32GM140911. The content is solely the authors’ responsibility and does not necessarily represent the official views of the National Institutes of Health. The performance evaluation was performed using the Wahab computing cluster at Old Dominion University. We would like to thank Min Dong for his assistance in creating the PODM docker image and
 %Angelos Angelopoulos for designing the docker container, Kevin Garner for providing usage guidelines for the ParaView software and 
 Emmanuel Billias for reviewing the manuscript and being a Beta-tester for the software.

\bibliographystyle{elsarticle-num}
\bibliography{references.bib}

\end{document}